\documentclass[12pt]{article}
\usepackage{mathrsfs}
\usepackage{epic,eepic,epsf,epsfig}
\usepackage{amsfonts,srcltx,mathrsfs}
\usepackage{amsmath}
\usepackage{amssymb}
\usepackage{amsbsy}
\usepackage{graphicx}
\usepackage{amsfonts}
\usepackage{color}

\newtheorem{lem}{Lemma}[section]
\newtheorem{theorem}[lem]{Theorem}

\newtheorem{cor}[lem]{Corollary}

\newtheorem{case}{Case}
\newtheorem{prop}[lem]{Proposition}

\def\a{\alpha} \def\b{\beta} \def\g{\gamma}

\def\di{\bigm|}  
\def\nd{\mathrel{\bigm|\kern-.7em/}}

\def\f{\noindent}

\def\PSL{\hbox{\rm PSL}}
\def\AGL{\hbox{\rm AGL}}
\def\SL{\hbox{\rm SL}}

\def\M{\hbox{\rm M}}

\def\S{\hbox{\rm S}}

\def\PSU{\hbox{\rm PSU}}
\def\Aut{\hbox{\rm Aut}}
\def\Inn{\hbox{\rm Inn}}
\def\Syl{\hbox{\rm Syl}}

\def\Aut{\hbox{\rm Aut}}\def\Out{\hbox{\rm Out}}
\def\N{\hbox{\rm N}}\def\C{\hbox{\rm C}}\def\O{\hbox{\rm O}}\def\S{\hbox{\rm S}}\def\A{\hbox{\rm A}}

\def\AGL{\hbox{\rm AGL}}
\def\PSp{\hbox{\rm PSp}}

\def\BiCay{\hbox{\rm BiCay}}\def\ASL{\hbox{\rm ASL}}\def\AGL{\hbox{\rm AGL}}

\def\M{\hbox{\rm M}}
\def\GL{\hbox{\rm GL}}
\def\SL{\hbox{\rm SL}}\def\O{\hbox{\rm O}}

\def\demo{\f {\bf Proof.}\hskip10pt}

\newcommand{\qed}{\mbox{\raisebox{0.7ex}{\fbox{}}} \vspace{4truemm}}
\def\mz{{\mathbb Z}}

\begin{document}

\title{Bipartite bi-Cayley graphs over metacyclic groups of odd prime-power order}

\author{ \\ Yi Wang, Yan-Quan Feng*\\
{\small\em Department of Mathematics, Beijing Jiaotong University, Beijing,
100044, P.R. China}\\}

\date{}
\maketitle

\footnotetext{*Corresponding author. E-mails: yiwang$@$bjtu.edu.cn, yqfeng$@$bjtu.edu.cn
}

\begin{abstract}
A graph $\Gamma$ is a bi-Cayley graph over a group $G$ if $G$ is a semiregular group of automorphisms of $\Gamma$ having two orbits. Let $G$ be a non-abelian metacyclic $p$-group for an odd prime $p$, and let $\Gamma$ be a connected bipartite bi-Cayley graph over the group $G$. In this paper, we prove that $G$ is normal in the full automorphism group $\Aut(\Gamma)$ of $\Gamma$ when $G$ is a Sylow $p$-subgroup of $\Aut(\Gamma)$. As an application, we classify  half-arc-transitive bipartite bi-Cayley graphs over the group $G$ of valency less than $2p$. Furthermore,  it is shown that there are no semisymmetric and no arc-transitive bipartite bi-Cayley graphs over the group $G$ of valency less than $p$.

\bigskip
\f {\bf Keywords:} bi-Cayley graph, half-arc-transitive graph, metacyclic group.
\end{abstract}

\section{Introduction}

All graphs considered in this paper are finite, connected, simple and undirected. For a graph $\Gamma$,
we use $V(\Gamma)$, $E(\Gamma)$, $A(\Gamma)$ and $\Aut(\Gamma)$ to denote its vertex set, edge set, arc set and full automorphism group, respectively. The graph $\Gamma$ is said to be {\em vertex-transitive}, {\em edge-transitive} or {\em arc-transitive} if $\Aut(\Gamma)$ acts transitively on $V(\Gamma)$, $E(\Gamma)$ or $A(\Gamma)$ respectively, {\em semisymmetric} if it is edge-transitive but not vertex-transitive, and {\em half-arc-transitive} if it is vertex-transitive, edge-transitive, but not arc-transitive.

Let $G$ be a permutation group on a set $\Omega$ and $\alpha \in \Omega$. Denote by $G_{\alpha}$ the stabilizer of $\alpha$ in $G$, that is, the subgroup of $G$ fixing the point $\alpha$. We say that $G$ is {\em semiregular} on $\Omega$ if $G_{\alpha}=1$ for every $\alpha \in \Omega$ and {\em regular} if $G$ is transitive and semiregular. A group $G$ is {\em metacyclic} if it has a normal subgroup $N$ such that both $N$ and  $G/N$ are cyclic.

Let $\Gamma$ be a graph with $G\leq \Aut(\Gamma)$. Then $\Gamma$ is said to be a {\em Cayley graph} over $G$ if $G$ is regular on $V(\Gamma)$ and a {\em bi-Cayley graph} over $G$ if $G$ is semiregular on $V(\Gamma)$ with two orbits. In particular, if $G$ is normal in $\Aut(\Gamma)$, the Cayley graph or bi-Cayley graph $\Gamma$ over $G$ is called {\em a normal Cayley graph} or {\em a normal bi-Cayley graph}, respectively.

It is well-known that Cayley graphs play an important role in the study of symmetry in graphs. However, graphs with various symmetries can be constructed by bi-Cayley graphs. For example, by using bi-Cayley graphs, several infinite families of semisymmetric graphs were constructed in~\cite{Du1, Du2, Lu}. Bi-Cayley graphs can also be used to construct non-Cayley vertex-transitive graphs, and the typical examples are the generalized Petersen graphs which are bi-Cayley graphs over cyclic groups. The smallest half-arc-transitive graph constructed in Bouwer~\cite{Bouwer} is also a bi-Cayley graph over a non-abelian metacyclic group of order $27$. In this paper, we construct a family of half-arc-transitive graphs by using bi-Cayley graphs.

In 1966, Tutte~\cite{Tutte} initiated an investigation of half-arc-transitive graphs by showing that a vertex- and edge-transitive graph with odd valency must be arc-transitive. A few years later, in order to answer Tutte's question of the existence of half-arc-transitive graphs of even valency, Bouwer~\cite{Bouwer} gave a construction of a $2k$-valent half-arc-transitive graph for every $k\geq 2$. One of the standard problems in the study of half-arc-transitive graphs is to classify such graphs of certain orders. Let $p$ be a prime. It is well known that there are no half-arc-transitive graphs of order $p$ or $p^2$, and by Cheng and Oxley~\cite{Cheng and Oxley}, there are no half-arc-transitive graphs of order $2p$. Alspach and Xu~\cite{Alspach} classified half-arc-transitive graphs of order $3p$ and Kutnar et al.~\cite{Kutnar} classified the half-arc-transitive graphs of order $4p$.
Despite all of these efforts, however, further classifications of half-arc-transitive graphs with general valencies seem to be very difficult,
and special attention has been paid to the study of half-arc-transitive graphs with small valencies, see~\cite{Conder1, Conder2, Feng-valency 4-2, Feng3-1, Marusic group action-2, Sparl, WF-1, XYWangF}. In fact, half-arc-transitive graphs have been extensively studied from different perspectives over decades by many authors; see, for example~\cite{Antoncic, Hujdurovic, Kutnar1, Li2, Potocnik, Spiga}.
Let $G$ be a non-abelian metacyclic $p$-group for an odd prime $p$. In this paper,
we classify  half-arc-transitive bipartite bi-Cayley graphs over the group $G$ of valency less than $2p$.

Our motivation comes partly from the work of Li and Sim~\cite{Li1,Li2}. Let $G$ be  a non-abelian metacyclic $p$-group for an odd prime $p$. In \cite{Li1}, the automorphism group of a Cayley graph over $G$ is characterized when $G$ is a Sylow $p$-subgroup of the Cayley graph, and by using this result, half-arc-transitive graphs over a metacyclic $p$-group of valency less than $2p$ were classified in~\cite{Li2}. In this paper, we first determine the automorphism group of a bi-Cayley graph over $G$ when $G$ is a Sylow $p$-subgroup of the bi-Cayley graph, and then using the result, we classified half-arc-transitive bipartite bi-Cayley graphs over $G$ of valency less than $2p$. Furthermore, we show that there is no semisymmetric and no arc-transitive bipartite bi-Cayley graphs over $G$ of valency less than $p$.

\section{Background results}\label{section 2}

Let $G$ be  a finite metacyclic $p$-group. Lindenberg~\cite{Linden} proved that the automorphism group of $G$ is a $p$-group when $G$ is nonsplit. The following proposition described the automorphism group of the remaining case when $G$ is split. It is easy to show that every non-abelian split metacyclic $p$-group $G$ for an odd prime $p$ has the following presentation:
$$ G_{\a,\b,\gamma}=\langle a,b\ |\ a^{p^\a}=1,\ b^{p^\b}=1,\ b^{-1}ab=a^{1+p^\gamma}\rangle,$$
\f where $\a,\b,\gamma$ are positive integers such that $0<\gamma<\a\leq \b+\gamma$. Let $n$ be a positive integer. Denote by $\mathbb{Z}_n$ the cyclic group of order $n$ as well as the ring of integers modulo $n$, and by $\mathbb{Z}_{n}^{*}$ the multiplicative group of the ring $\mathbb{Z}_n$ consisting of numbers coprime to $n$.

\begin{prop}\label{prop=$p-1$-subgroup}{\rm(\cite[Theorem~2.8]{Li2})}
For an odd prime $p$, we have that
$$|\Aut(G_{\a,\b,\gamma})|=(p-1)p^{{\rm min}(\a,\b)+{\rm min}(\b,\gamma)+\b+\gamma-1}.$$
Moreover, all Hall $p^{'}$-subgroups of $\Aut(G_{\a,\b,\gamma})$ are conjugate and isomorphic to $\mathbb{Z}_{p-1}$. In particular, the map $\theta: a\mapsto a^\varepsilon,\ b\mapsto b,$ induces an automorphism of $G_{\a,\b,\gamma}$ of order $p-1$,  where $\varepsilon$ is an element of order $p-1$ in $\mz_{p^\a}^*$.
\end{prop}

A $p$-group $G$ is said to be {\em regular} if for any $x,y\in G$ there exists $d_i\in \langle x,y\rangle^{'}$ such that $x^{p}y^{p}=(xy)^{p}\Pi d_{i}^{p}$. If $G$ is metacyclic, then the derived group $G^{'}$ is cyclic, and hence $G$ is regular by \cite[Kapitel III, 10.2 Satz]{Huppert}. For regular $p$-groups, the following proposition holds by \cite[Kapitel III, 10.8 Satz]{Huppert}.

\begin{prop}\label{prop=element order}
Let $G$ be a metacyclic $p$-group for an odd prime $p$. If $|G^{'}|=p^{n}$, then for any $m\geq n$, we have
$$(xy)^{p^m}=x^{p^m}y^{p^m},\ \forall x,y\in G.$$
\end{prop}

It is easy to check that $|G_{\a,\b,\gamma}^{'}|=p^{\a-\gamma}$. Let $b^{m}a^{n}\in G_{\a,\b,\gamma}$ with $m\in \mz_{p^\b}$ and $n\in \mz_{p^\a}$, and denote by $o(b^ma^n)$ the order of $b^ma^n$. Since $\a-\gamma\leq \b$, Proposition~\ref{prop=element order} implies that if $(p,m)=1$ then $o(b^{m}a^{n})=\max\{o(a^n),p^\b\}$, and if $\b<\a$ and $p\di n$, then $o(b^{m}a^{n})\leq p^{\a-1}$, which will be used later.
Let $G$ be a finite group. Denote by $N\leq G$ if $N$ is a subgroup of $G$, and by $N<G$ if $N$ is a proper subgroup of $G$. The following proposition determines non-abelian simple groups having a proper subgroup of index prime-power order.

\begin{prop}{\rm(\cite[Theorem 1]{Guralnick})} \label{prop=subgroups in a simple group}
Let $T$ be a non-abelian simple group with $H<T$ and $|T:H|=p^a$, p prime. Then one of the following holds.
\begin{enumerate}
\item[$1$.] $T=\PSL(n,q)$ and $H$ is the stabilizer of a line or hyperplane. Furthermore, $|T:H|=(q^n-1)/(q-1)=p^a$ and $n$ must be a prime.
\item[$2$.]  $T=\A_n$ and $H\cong \A_{n-1}$ with $n=p^a$.
\item[$3$.] $T=\PSL(2,11)$ and $H\cong \A_5$.
\item[$4$.] $T=\M_{23}$ and $H\cong \M_{22}$ or $T=\M_{11}$ and $H\cong \M_{10}$.
\item[$5$.] $T=\PSU(4,2)\cong \PSp(4,3)$ and $H$ is the parabolic subgroup of index $27$.
\end{enumerate}
\end{prop}

The following result is an immediate consequence of Corollary~2 in Guralnick~\cite{Guralnick}.

\begin{prop}  \label{prop=2-transitive}
Let $T$ be a non-abelian simple group acting transitively on $\Omega$ with $p^{l}$ elements for a prime $p$. If $p$ does not divide the order of a point-stabilizer in $T$, then $T$ acts $2$-transitively on $\Omega$.
\end{prop}

It is well-know that $\GL(d,q)$ has a cyclic group of order $q^d-1$, the so called Singer-Zyklus, which also induces a cyclic group on $\PSL(d,q)$.

\begin{prop}{\rm(\cite[Kapitel II, 7.3 Satz]{Huppert})}
\label{prop=cyclic subgroup}
The group $G=\GL(d,q)$ contains a cyclic subgroup of order $q^d-1$, and it induces a cyclic subgroup of order $\frac{q^d-1}{(q-1)(q-1,d)}$ on $\PSL(d,q)$.
\end{prop}

Let $G$ and $E$ be two groups. We call an extension $E$ of $G$ by $N$ a {\em central extension} of $G$ if $E$ has a central subgroup $N$ such that $E/N\cong G$,
and if further $E$ is perfect, that is, the derived group $E'=E$, we call $E$ a {\em covering group} of $G$. Schur~\cite{Schur} proved that for every non-abelian simple group $G$ there is a unique maximal covering group $M$ such that every covering group of $G$ is a factor
group of $M$ (also see \cite[Chapter 5, Section 23]{Huppert}). This group $M$ is called the {\em full covering group} of $G$, and the center of $M$ is the {\em Schur multiplier} of $G$, denoted by $\M(G)$. For a group $G$, we denote by $\Out(G)$ the outer automorphism group of $G$, that is, $\Out(G)=\Aut(G)/\Inn(G)$, where $\Inn(G)$ is the inner automorphism group of $G$ induced by conjugation.

The following proposition is about outer automorphism group and Schur multiplier of a non-abelian simple group with a proper subgroup of prime-power index.

\begin{prop} {\rm(\cite[Lemma~2.3]{Li1})} \label{prop=outautomorphism}
Let $p$ be an odd prime and let $T$ be a non-abelian simple group which has a subgroup $H$ of index $p^{l}>1$. Then

\begin{enumerate}
\item[$1$.]$p \nmid|\M(T)|$;
\item[$2$.]either $p \nmid|\Out(T)|$ or $T\cong\PSL(2,8)$ and $p^{l}=3^{2}$.
\end{enumerate}
\end{prop}

For a group $G$ and a prime $p$, denote by $\O_p(G)$ the largest normal $p$-subgroup of $G$, and by $\O_{p'}(G)$ the maximal normal subgroup of $G$ whose order is not divisible by $p$. The next proposition is about transitive permutation groups of prime-power degree.

\begin{prop} {\rm(\cite[Lemma~2.5]{Li1})}\label{prop=propersubgroup}
Let $p$ be a prime, and let $A$ be a transitive permutation group of $p$-power degree. Let $B$ be a nontrivial subnormal subgroup of $A$. Then $B$ has a proper subgroup of $p$-power index, and $\O_{p'}(B)=1$. In particular, $\O_{p'}(A)=1$.
\end{prop}

A group $G$ is said to be {\em a central product} of
the subgroups $H_{1},\ldots, H_{n}$ ($n\geq 2$) of $G$ if $G=H_{1}\cdots H_{n}$ and for any $i\neq j$, $H_{i}$ and $H_{j}$ commute elementwise. A group $G$ is called {\em quasisimple} if $G^{'}=G$ and $G/Z(G)$ is a non-abelian simple group, where $Z(G)$ is the centralizer of $G$, and is called {\em semisimple} if $G^{'}=G$ and $G/Z(G)$ is a direct product of non-abelian simple groups. Clearly, a quasisimple group is semisimple, and the identity group is semisimple, but not quasisimple.

\begin{prop}{\rm(\cite[Theorem~6.4]{Suzuki2})} \label{prop=semisimplle}
A central product of two semisimple groups is also semisimple. Any semisimple group can be decomposed into a central product of quasisimple groups, and this set of quasisimple groups is uniquely determined.
\end{prop}

A subnormal quasisimple subgroup of a group $G$ is called a {\em component} of $G$. By \cite[6.9(iv), p.~450]{Suzuki2}, any two distinct components of $G$ commute elementwise, and by Proposition~\ref{prop=semisimplle}, the product of all components of $G$ is semisimple, denoted by $E(G)$, which is characteristic in $G$.
We use $F(G)$ to denote the {\em Fitting subgroup} of $G$, that is, $F(G)= \O_{p_{1}}(G)\times \O_{p_{2}}(G)\times\cdots\times \O_{p_{t}}(G)$, where $p_{1}, p_{2},\ldots, p_{t}$ are all distinct prime factors of $|G|$. Set $F^{*}(G)=F(G)E(G)$ and call $F^{*}(G)$ the {\em generalized Fitting subgroup} of $G$. The following is one of the most significant properties of $F^{*}(G)$.

\begin{prop}{\rm(\cite[Theorem~6.11]{Suzuki2})} \label{prop=GeneralizedFitting}
For any finite group $G$, we have $$\C_{G}(F^{*}(G))\leq F^{*}(G).$$
\end{prop}

For a group $G$ and a subgroup $H$ of $G$, let $\C_G(H)$ be the centralizer of $H$ in $G$. An {\it action} of a group $G$ on a set $\Omega$ is a homomorphism from $G$ to the symmetric group $S_{\Omega}$ on $\Omega$. We denote by $\Phi(G)$ the {\em Frattini subgroup} of $G$, that is, the intersection of all maximal subgroups of $G$. Note that for a prime $p$, $\O_p(G)$ is a $p$-group and $\O_p(G)/\Phi(\O_p(G))$ is an elementary abelian $p$-group. Thus, $\O_p(G)/\Phi(\O_p(G))$ can be viewed as a vector space over the field $\mz_p$.
The following lemma considers a natural action of a group $G$ on the vector space $\O_p(G)/\Phi(\O_p(G))$.

\begin{prop}{\rm(\cite[Lemma~2.9]{Wang})}
\label{prop=kernel}
For a finite group $G$ and a prime $p$, let $H=\O_{p}(G)$ and $V=H/\Phi(H)$. Then $G$ has a natural action on $V$, induced by conjugation via elements of $G$ on $H$. If $\C_{G}(H)\leq H$, then $H$ is the kernel of this action of $G$ on $V$.
\end{prop}

Let $T$ act on two sets $\Omega$ and $\Sigma$, and these two actions are {\em equivalent} if there is a bijection $\lambda : \Omega \mapsto \Sigma$ such that $$(\alpha^{t})^\lambda=(\alpha^\lambda)^t \ {\rm for\ all}\ \alpha \in \Omega \ {\rm and}\ t\in T.$$ When the two actions are transitive there is a simple criterion for deciding whether or not they are equivalent.

\begin{prop}{\rm(\cite[Lemma~1.6B]{Dixon})}
\label{prop=equivalent action}
Assume that the group $T$ acts transitively on the two sets $\Omega$ and $\Sigma$, and let $W$ be a stabilizer of a point in the first action. Then the actions are equivalent if and only if $W$ is the stabilizer of some point in the second action.
\end{prop}

For a group $G$ and two subgroups $H$ and $K$ of $G$, we consider the actions of $G$ on the right cosets of $H$ and $K$ by right multiplication. The stabilizers of $Hx$ and $Ky$ are $H^x$ and $K^y$, respectively. By Proposition~\ref{prop=equivalent action},  these two right multiplication actions are equivalent if and only if $H$ and $K$ are conjugate in $G$.

\section{Automorphism groups}\label{section 3}

Let $p$ be an odd prime and $G$ a non-abelian metacyclic $p$-group for an odd prime $p$. If $\Gamma$ is a connected Cayley graph over a group $G$ and $G\unlhd \Aut(\Gamma)$, then $\Aut(\Gamma)$ is known by Godsil~\cite{Godsil}. If $\Gamma$ is a connected bi-Cayley graph over $G$ and $G\unlhd \Aut(\Gamma)$, then $\Aut(\Gamma)$ is known by Zhou and Feng~\cite{Zhou1} (also see Proposition~\ref{prop=Normalizer}). The main purpose of this section is to determine automorphism group of a connected  bipartite bi-Cayley graph over $G$ by proving that $G$ is normal in the full automorphism group of the graph.

Let $\Gamma_N$ be the {\em quotient graph} of a graph $\Gamma$ with respect to $N\leq \Aut(\Gamma)$, that is, the graph having the orbits of $N$ as vertices with two orbits $O_1$, $O_2$ adjacent in $\Gamma_{N}$ if there exist some $u\in O_1$ and $v\in O_2$ such that $\{u,v\}$ is an edge in $\Gamma$. Denote by $[O_1]$ the induced subgraph of $\Gamma$ by $O_1$, and by $[O_1,O_2]$ the spanning subgraph of  $[O_1\cup O_2]$ with edge set $\{\{u,v\}\in E(\Gamma)\ |\ u\in O_1, v\in O_2\}$.
The following is the main result of this section.

\begin{theorem}\label{maintheorem}
Let $G$ be a non-abelian metacyclic group of order $p^{n}$ for an odd prime $p$ and a positive integer $n$, and let $\Gamma$ be a connected bipartite bi-Cayley graph over $G$. Assume that $G$ is a Sylow $p$-subgroup of $\Aut(\Gamma)$. Then $G\unlhd \Aut(\Gamma)$.
\end{theorem}

\demo  Let $A=\Aut(\Gamma)$ and let $W_0$ and $W_1$ be the two parts of the bipartite graph $\Gamma$. Then $\{W_0, W_1\}$ is a complete block system of $A$ on $V(\Gamma)$. Let $A^*$ be the kernel of $A$ on $\{W_0, W_1\}$, that is, the subgroup of $A$ fixing $W_0$ and $W_1$ setwise. It follows $A^*\unlhd A$, $A/A^*\leq \mz_2$ and $\Syl_p(A)=\Syl_p(A^*)$ because $p$ is odd, where $\Syl_p(A)$ and $\Syl_p(A^*)$ are the sets of Sylow $p$-subgroups of $A$ and $A^*$ respectively. Since $G$ is a Sylow $p$-subgroup of $A$, we have $G\in \Syl_p(A^*)$. In particular, $p^n\mid |A|$ and $p^{n+1}\nmid |A|$, that is,  $p^n\parallel |A|$. Since $G\leq A^*$, the two orbits of $G$ are exactly $W_0$ and $W_1$, and $G$ is regular on both $W_0$ and $W_1$. By Frattini argument~\cite[Kapitel I, 7.8 Satz]{Huppert}, $A^*=GA^*_u$ for any $u\in V(\Gamma)$, implying that $A^*_u$ is a $p'$-group. Clearly, $A_u=A^*_u$ is also a $p'$-group. Since Sylow $p$-subgroups of $A$ are conjugate, every $p$-subgroup of $A$ is semiregular on both $W_0$ and $W_1$.

If the kernel of $A^*$ on $W_0$ (resp. $W_1$) is unfaithful, it has order  divisible by $p$ as $|W_1|=p^n$, and so $p^{n+1}\mid |A^*|$, a contradiction. Thus, $A^*$ acts faithfully on $W_0$ (resp. $W_1$).

\medskip
\f{\bf Claim~1:} Any minimal normal subgroup $N$ of $A^*$ is abelian.

Suppose that $N$ is non-abelian. Then $N\cong T_{1}\times \cdots \times T_{k}$ with $k\geq 1$, where $T_{i}\cong T$ is a non-abelian simple group. By Proposition~\ref{prop=propersubgroup}, $|N|$ is divisible by $p$. Since $G \in \Syl_{p}(A^*)$, we have $G \cap N \in \Syl_p(N)$, and hence $G \cap N=P_{1}\times \cdots\times P_{k}$ for some $P_{i} \in \Syl_p(T_{i})$. Then $P_i \neq 1$ for $1 \leq i\leq k$, as otherwise $p\nmid |N|$. Since $G$ is metacyclic and $G\cap N\unlhd G$, $G\cap N$ is metacyclic, forcing $k\leq 2$.

Set $\Omega=\{T_{1},\ldots, T_{k}\}$ and write $B=\N_{A^*}(T_{1})$. By considering the conjugation action of $A^*$ on $\Omega$, we have $B\unlhd A^*$ as $k\leq 2$, and hence $A^*/B\lesssim \S_2$. Thus, each Sylow $p$-subgroup of $B$ is also a Sylow $p$-subgroup of $A^*$, implying that $B$ is transitive on both $W_0$ and $W_1$.

Let $\Gamma_{T_{1}}$ be the quotient graph of $\Gamma$ with respect to $T_{1}$. Since $T_1\unlhd B$, all orbits of $T_1$ on $W_0$ have the same length, and the length must be a $p$-power as $|W_0|=p^n$, so that it is the order of a Sylow $p$-subgroup of $T_1$ because each $p$-subgroup is semiregular.  Similarly, all orbits of $T_1$ on $W_1$ have the same length and it is also the order of a Sylow $p$-subgroup of $T_1$. Thus, $V(\Gamma_{T_{1}})=\{\Delta_{1},\ldots, \Delta_{s}, \Delta_{1}^{'},\ldots, \Delta_{s}^{'}\}$, the set of all orbits of $T_1$, with $W_0=\Delta_1\cup \cdots \cup \Delta_s$ and $W_1=\Delta_{1}^{'}\cup \cdots \cup \Delta_{s}^{'}$. Furthermore, for any $1\leq i,j\leq s$ we have $|\Delta_{i}|=|\Delta_{j}^{'}|= p^{m}$ for some $1\leq m\leq n$, and hence $s=p^{n-m}$. Since $T_1\unlhd B$, $B$ has a natural action on $V(\Gamma_{T_1})$ and let $K$ be the kernel of this action. Clearly, $T_{1}\leq K$.  Recall that $p\nmid |A_u|$ for any $u\in V(\Gamma)$. Then $p\nmid |(T_1)_u|$, and by Proposition~\ref{prop=2-transitive}, $T_{1}$ is $2$-transitive on each $\Delta_i$ or $\Delta_i'$. Since $p^n \parallel |A^*|$, we have $p^n \parallel |B|$, implying that $p^m\parallel |K|$. Since $(T_1)_u$  is a proper subgroup of $T_1$ of index $p$-power, Proposition~\ref{prop=outautomorphism} implies that either $T_1=\PSL(2,8)$ with $p^m=3^2$, or $p\nmid|\Out(T_1)|$. To finish the proof of the claim, we only need to prove that these two cases are impossible.

\medskip
\f{\bf Case~1:} $T_1=\PSL(2,8)$ with $p^m=3^2$.

In this case, $|\Delta_i|=|\Delta_{j}^{'}|=9$. If $s=1$ then $|G|=p^m=3^2$, contrary to the fact that $G$ is non-abelian. Now assume $s\geq 2$. By Atlas~\cite{Conway},  $\PSL(2,8)$ has only one conjugate class of subgroups of index $9$, and by Proposition~\ref{prop=equivalent action}, $T_1$ acts equivalently on $\Delta_i$ and $\Delta_{j}^{'}$.

\begin{figure}[h]
\begin{center}
\includegraphics[height=6.2cm,width=15cm]{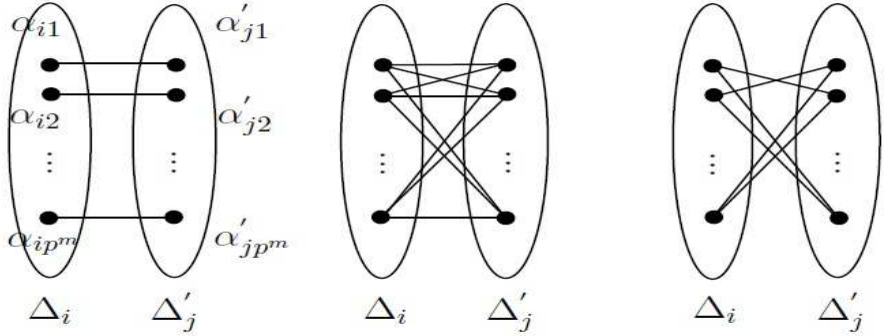}
\end{center}
\vskip -1 cm
\caption{The subgraph $[\Delta_i,\Delta_{j}^{'}]$}
\end{figure}

Set $\Delta_i=\{\a_{i1},\a_{i2},\ldots,\a_{i9}\}$ and $\Delta_{j}^{'}=\{\a_{j1}^{'},\a_{j2}^{'},\ldots,\a_{j9}^{'}\}$ for $1\leq i,j\leq 3^{n-2}$.  Note that $T_1$ is $2$-transitive on $\Delta_i$ and  $\Delta_{j}^{'}$.
Since $T_1$ acts equivalently on $\Delta_i$ and $\Delta_{j}^{'}$, by Proposition~\ref{prop=equivalent action}, we may assume that $(T_1)_{\a_{il}}=(T_1)_{\a_{jl}^{'}}$ for any $1\leq i,j\leq 3^{n-2}$ and $1\leq l\leq 3^2$. The subgraph $[\Delta_i,\Delta_{j}^{'}]$ is either a null graph, or one of the three graphs in Figure~1 because $(T_1)_{\a_{il}}=(T_1)_{\a_{jl}^{'}}$ acts transitively on both $\Delta_i\backslash \{\a_{il}\}$ and $\Delta_{j}^{'}\backslash \{\a_{jl}^{'}\}$. The three graphs are of edge sets $\{\{\alpha_{il}, \alpha_{jl}^{'}\}\ |\ 1\leq l\leq 3^2\}$, $\{\{\alpha_{ik}, \alpha_{jl}^{'}\}\ |\ 1\leq k, l\leq 3^2\}$ and $\{\{\alpha_{ik}, \alpha_{jl}^{'}\}\ |\ 1\leq k, l\leq 3^2, k\not=l\}$.

For any $g\in \S_{9}$, define a permutation $\sigma_g$ on $V(\Gamma)$ by $(\alpha_{il})^{\sigma_g}=\alpha_{il^g}$ and $(\a_{jl}^{'})^{\sigma_g}=\a_{jl^g}^{'}$ for any $1\leq i,j\leq 3^{n-2}$ and $1\leq l \leq 3^2$. Then $\sigma_g$ fixes each $\Delta_i$ and $\Delta_{j}^{'}$, and permutes the elements of $\Delta_i$ and $\Delta_{j}^{'}$ in the `same way' for each $1\leq i,j\leq 3^{n-2}$. Since $[\Delta_i, \Delta_{j}^{'}]$ is either a null graph, or one graph in Figure 1, $\sigma_g$ induces an automorphism of $[\Delta_i, \Delta_{j}^{'}]$, for all $1\leq i, j\leq 3^{n-2}$. Also $\sigma_g$ induces automorphisms of $[\Delta_i]$ and $[\Delta_{j}^{'}]$ for all $1\leq i,j\leq 3^{n-2}$ because  $[\Delta_i]$ and  $[\Delta_{j}^{'}]$ have no edge ($\Gamma$ is bipartite). It follows that $\sigma_g \in \Aut(\Gamma)$. Thus, $L:=\{ \sigma_g\ |\ g\in \S_{9} \} \leq \Aut(\Gamma)$ and $L \cong \S_{9}$.
Recall that $K$ is the kernel of $B$ acting on $V(\Gamma_{T_1})$.
Since $L$ fixes each $\Delta_i$ and $\Delta_{j}^{'}$, we have $L\leq K$, and since $3^3\di |L|$, we have $3^3\di |K|$, contrary to the fact that $3^2\parallel |K|$.

\medskip
\f{\bf Case~2:} $p\nmid |\Out(T_1)|$.

Since $B/T_1 \C_{B}(T_1)\lesssim \Out(T_1)$, we have $p^n\parallel |T_1 \C_{B}(T_1)|$. Since $T_1$ is non-abelian simple, $T_1\cap \C_{B}(T_1)=1$ and hence $T_1\C_{B}(T_1)=T_1\times \C_{B}(T_1)$. If $p\di |\C_{B}(T_1)|$, then $G$ is conjugate to $Q_1\times Q_2$, where $Q_1\in \Syl_p(T_1)$ and $Q_2\in \Syl_p(\C_{B}(T_1))$. Since $G$ is metacyclic, $G$ can be generated by two elements, and since $G$ is a $p$-group, any minimal generating set of $G$ has cardinality $2$. It follows that both $Q_1$ and $Q_2$ are cyclic, and so $G$ is abelian, a contradiction. Thus, $p\nmid |\C_{B}(T_1)|$ and hence $p^n\parallel |T_1|$, forcing $s=1$. Furthermore, $W_0=\Delta_1$, $W_1=\Delta_{1}^{'}$ and $T_1$ is $2$-transitive on both $W_0$ and $W_1$. Note that $(T_1)_u$ is proper subgroup of $T_1$ of index $p^n$. Since $G$ is a Sylow $p$-subgroup of $A$ of order $p^n$, all Sylow $p$-subgroups of $T_1$ are also Sylow $p$-subgroups of $A$, and so they are isomorphic to $G$. Without loss of generality, we may assume $G\leq T_1$. By Proposition~\ref{prop=subgroups in a simple group}, $T_1=\PSL(2,11)$, $\M_{11}$, $\M_{23}$, $\PSU(4,2)$, $\A_{p^n}$, or $\PSL(d,q)$ with $\frac{q^d-1}{q-1}=p^n$ and $d$ a prime.

Suppose $T_1=\PSL(2,11)$, $\M_{11}$ or $\M_{23}$. By Proposition~\ref{prop=subgroups in a simple group}, $|W_0|=|W_1|=11$, $11$, or $23$ respectively, and hence $|G|=11$, $11$ or $23$, contrary to the fact that $G$ is non-abelian.

Suppose $T_1=\PSU(4,2)$ or $\A_{p^n}$. For the former, $T_1$ has one conjugate class of subgroups of index $27$ by Atlas~\cite{Conway}, and for the latter,  $T_1$ has  one conjugate class of subgroups of index $p^n$. By Proposition~\ref{prop=equivalent action}, $T_1$ acts equivalently on $W_0$ and $W_1$, and since $\Gamma$ is connected, the $2$-transitivity of $T_1$ on $W_0$ and $W_1$ implies that $\Gamma\cong K_{p^n,p^n}$ or $K_{p^n,p^n}-p^n K_2$. Then $A=\S_{p^n}\wr \S_2$ or $\S_{p^n}\times \mathbb{Z}_2$ respectively. Since $G$ is non-abelian, we have $n\geq 3$, and so $p^{n+1}\di |A|$, a contradiction.

Suppose $T_1=\PSL(d,q)$ with $\frac{q^d-1}{q-1}=p^n$ and $d$ a prime. By Proposition~\ref{prop=cyclic subgroup}, $T_1$ has a cyclic subgroup of order $\frac{q^d-1}{(q-1)(q-1,d)}$. Since $d$ is a prime, either $(q-1,d)=1$ or $(q-1,d)=d$. Note that $(q-1,d)\di \frac{q^d-1}{q-1}$. If $(q-1,d)=d$ then $d=p$ and $p\di (q-1)$. Since $p\geq 3$ and $p^{2}\di (q^2-1)(q-1)$, we have
$p^{n+1}\di \frac{(q^p-1)(q^p-q)\cdots (q^p-q^{p-1})}{(q-1)(q,d)}$, that is, $p^{n+1}\di |T_1|$, a contradiction.
If $(q-1,d)=1$ then $T_1$ has a cyclic subgroup of order $\frac{q^d-1}{q-1}=p^n$, contrary to the fact that $G$ is non-abelian. This completes the proof of Claim~1.

\medskip
\f{\bf Claim~2:} $\C_{A^*}(\O_{p}(A^*))\leq \O_{p}(A^*)$.

Suppose that $D$ is a component of $A^*$, that is, a subnormal quasisimple subgroup of $A^*$. Then $D=D'$ and $D/Z(D) \cong T$, a non-abelian simple group. By Proposition~\ref{prop=propersubgroup}, $D$ has a proper subgroup $C$ of $p$-power index and $Z(D)$ is a $p$-group. Since $|D:C|=|D:CZ(D)|\cdot |CZ(D):C|$, we have that $|D: CZ(D)|$ is a $p$-power. If $D=CZ(D)$ then $D=D^{'}=C^{'}$, contrary to the fact that $C$ is a proper subgroup of $D$. Thus, $CZ(D)\neq D$. Since $|D/Z(D): CZ(D)/Z(D)|=|D: CZ(D)|$, we have that $D/Z(D)$ has a proper subgroup $CZ(D)/Z(D)$ of $p$-power index. By Proposition~\ref{prop=outautomorphism}(1), $p\nmid |M(D/Z(D))|$ and hence $p\nmid |Z(D)|$. Since $Z(D)$ is a $p$-group, $Z(D)=1$ and $D \cong T$. Recall that $E(A^*)$ is the product of all components of $A^*$. Then $D\leq E(A^*)$ and since $D\cong T$, $D$ is a direct factor of $E(A^*)$. Clearly, $D^a$ is also a direct factor of $E(A^*)$ for any $a\in A^*$. It follows that $A^*$ contains a minimal normal subgroup which is isomorphic to $T^\ell$ with $\ell\geq 1$, contrary to Claim~1. Thus, $A^*$ has no component and $E(A^*)=1$. It follows that the generalized Fitting subgroup $F^*(A^*)=F(A^*)$. By Proposition~\ref{prop=propersubgroup}, $\O_{p^{'}}(A^*)=1$ and hence $F^{*}(A^*)=\O_{p}(A^*)$. By Proposition~\ref{prop=GeneralizedFitting}, $\C_{A^*}(\O_{p}(A^*))\leq \O_{p}(A^*)$, as claimed.

\medskip

Now we are ready to finish the proof. To do it, it suffices to show $G\unlhd A^*$ as  $|A:A^*|\leq 2$.

Let $H=\O_p(A^*)$. By Claim~1, $H\not=1$. Write $\overline{H}=H/\Phi(H)$ and $\overline{A^*}=A^*/\Phi(H)$. Then $\O_p(A^*/H)=1$, and since $G\in \Syl_p(A^{*})$, $H\leq G$. By Claim~2 and Proposition~\ref{prop=kernel}, $A^*/H\leq \Aut(\overline{H})$. Since $G$ is metacyclic, $\overline{H}=\mathbb{Z}_{p}$ or $\mathbb{Z}_{p}\times \mathbb{Z}_{p}$.

Assume $\overline{H}= \mathbb{Z}_{p}$. Then $A^*/H\leq \mathbb{Z}_{p-1}$, and  $G=H\unlhd A^*$, as required.

Assume $\overline{H}=\mathbb{Z}_{p}\times \mathbb{Z}_{p}$. Then $A^*/H\leq \GL(2,p)$. If $p\nmid |A^*/H|$ then $G=H\unlhd A^*$, as required. Suppose $p\di |A^*/H|$. We finish the proof by showing that this is impossible.

Since $p\parallel |\GL(2,p)|$, we have $p\parallel|A^*/H|$, and since $\Syl_p(\SL(2,p))= \Syl_p(\GL(2,p))$, we have $\Syl_p(A^*/H)\subseteq \Syl_p(\SL(2,p))$.
Note that $A^*/H\cdot \SL(2,p)\leq \GL(2,p)$. Then $p\parallel|A^*/H\cdot \SL(2,p)|$ and $p\di |(A^*/H)\cap \SL(2,p)|$. Since $\O_p(A^*/H)=1$, $A^*/H$ has at least two Sylow $p$-subgroups, and hence $(A^*/H)\cap \SL(2,p)$ has at least two Sylow $p$-subgroups, implying that $(A^*/H)\cap \SL(2,p)$ has no normal Sylow $p$-subgroups. By {\rm(\cite[Theorem~6.17]{Suzuki1})}, $(A^*/H)\cap\SL(2,p)$ contains $\SL(2,p)$, that is, $\SL(2,p)\leq A^*/H\leq \GL(2,p)$.
In particular, the induced faithful representation of $A^*/H$ on the linear space $\overline{H}$ is irreducible, and hence  $\overline{H}$ is a minimal normal subgroup of  $\overline{A^*}$.

Since $\overline{A^*}\big/ \overline{H}\cong A^*/H$, we have $\SL(2,p)\leq\overline{A^*}/\overline{H}\leq \GL(2,p)$. Write $\overline{R}/\overline{H}=Z(\overline{A^*}\big/ \overline{H})$. Then $\overline{R}\unlhd \overline{A^*}$ and $1\not=\overline{R}/\overline{H}$ is a $p'$-group. Since $\overline{H}\unlhd \overline{R}$ and $\overline{H}\in\Syl_p(\overline{R})$, the Schur-Zassenhaus Theorem~\cite[Theorem~8.10]{Suzuki1} implies that there is a $p'$-group $\overline{V}\leq \overline{R}$ such that $\overline{R}=\overline{H}\overline{V}$ and all Hall $p^{'}$-subgroup of $\overline{R}$ are conjugate. Note that $\overline{V}\not=1$.  By Frattini argument~\cite[Kapitel I, 7.8 Satz]{Huppert}, $\overline{A^*}=\overline{R} N_{\overline{A^*}}(\overline{V})=\overline{H}N_{\overline{A^*}}(\overline{V})$.
Since $\overline{H}$ is abelian, $\overline{H}\cap N_{\overline{A^*}}(\overline{V})\unlhd \overline{A^*}$, and by the minimality of $\overline{H}$, $\overline{H}\cap N_{\overline{A^*}}(\overline{V})=\overline{H}$ or $1$. If $\overline{H}\cap N_{\overline{A^*}}(\overline{V})=\overline{H}$
then $\overline{H}\leq N_{\overline{A^*}}(\overline{V})$ and $\overline{A^*}=\overline{H}N_{\overline{A^*}}(\overline{V})=N_{\overline{A^*}}(\overline{V})$, that is, $\overline{V}\unlhd \overline{A^*}$. This implies that $\O_{p'}(\overline{A^*})\not=1$, contrary to Proposition~\ref{prop=propersubgroup}.
If $\overline{H}\cap N_{\overline{A^*}}(\overline{V})=1$ then  $\overline{A^*}=\overline{H}N_{\overline{A^*}}(\overline{V})$ implies  $\ASL(2,p) \leq \overline{A^*}\leq \AGL(2,p)$ as $\SL(2,p)\leq\overline{A^*}/\overline{H}\leq \GL(2,p)$. It follows that a Sylow $p$-subgroup of $\overline{A^*}$ is not metacyclic. On the other hand, since both normal subgroups and quotient groups of a metacyclic group are metacyclic, a Sylow $p$-subgroups of  $\overline{A^*}$ is metacyclic because a Sylow $p$-subgroup of $A^*$ is metacyclic, a contradiction. \hfill\qed

\section{Edge-transitive bipartite bi-Cayley graphs}\label{section 4}

A connected edge-transitive graph is semisymmetric, arc-transitive, or half-arc-transitive.
As an application of Theorem~\ref{maintheorem}, we prove that
there is no connected semisymmetric or arc-transitive bipartite bi-Cayley graph over a non-abelian metacyclic $p$-group with valency less than $p$. Moreover, we classify connected half-arc-transitive bipartite bi-Cayley graphs over a non-abelian metacyclic $p$-group with valency less than $2p$.

Let $G$ be a group and let $R$, $L$ and $S$ be subsets of $G$ such that $R=R^{-1}$, $L=L^{-1}$, $R\cup L$ does not contain the identity element of $G$ and $S$ contains the identity element of $G$. Define the graph $\BiCay(G,R,L,S)$ to have vertex set the union of the {\em right part} $W_0=\{g_0\ |\ g\in G\}$ and the {\em left part} $W_1=\{g_1\ |\ g\in G\}$, and edge set the union of the {\em right edges} $\{\{h_0,g_0\}\ |\ gh^{-1}\in R\}$, the {\em left edges} $\{\{h_1,g_1\}\ |\ gh^{-1}\in L\}$ and the {\em spokes} $\{\{h_0,g_1\}\ |\ gh^{-1}\in S\}$. The graph $\BiCay(G,R,L,S)$ is connected if and only if $G=\langle R\cup L\cup S\rangle$, and  $\BiCay(G,R,L,S)\cong \BiCay(G,R^\theta,L^\theta,S^\theta)$ for any $\theta\in \Aut(G)$.

On the other hand, let $\Gamma$ be a bi-Cayley graph over $G$, where $G$ is a semiregular group of automorphisms of $\Gamma$ with two orbits. Then $\Gamma$ can be realized as above, that is,  $\Gamma\cong \BiCay(G,R,L,S)$ for some subsets $R$, $L$ and $S$ of $G$ satisfying $R=R^{-1}$, $L=L^{-1}$,  $R\cup L$ does not contain the identity element of $G$ and $S$ contains the identity element of $G$. Thus, $\BiCay(G,R,L,S)$ is also called a {\em bi-Cayley graph} over $G$ relative to $R$, $L$ and $S$.

Let $\Gamma=\BiCay(G,R,L,S)$ be a connected bi-Cayley graph over a group $G$. For $g\in G$, define a permutation $\hat{g}$ on $V(\Gamma)=W_0 \cup W_1$ by the rule
$$h_{i}^{\hat{g}}=(hg)_i,\ \forall i\in \mathbb{Z}_2,\ h,g\in G.$$
It is easy to check that $\hat{g}$ is an automorphism of $\Gamma$ and  $\hat{G}=\{\hat{g}\ |\ g\in G\}$ is a semiregular group of automorphisms of $\Gamma$ with two orbits.

For an automorphism $\theta$ of $G$ and $x,y,g\in G$, define two permutations on $V(\Gamma)=W_0 \cup W_1$ as following:
$$\delta_{\theta,x,y}:\ h_0\mapsto (xh^\theta)_1,\ h_1\mapsto (yh^\theta)_0,\ \forall h\in G,$$
$$\sigma_{\theta,g}:\ h_0\mapsto (h^\theta)_0,\ h_1\mapsto (gh^\theta)_1,\ \forall h\in G.$$
Set
$$I:=\{\delta_{\theta,x,y}\ |\ \theta\in \Aut(G)\ s.t.\ R^\theta=x^{-1}Lx,\ L^\theta=y^{-1}Ry,\ S^\theta=y^{-1}S^{-1}x\},$$
$$F:=\{\sigma_{\theta,g}\ |\ \theta\in \Aut(G)\ s.t.\ R^\theta=R,\ L^\theta=g^{-1}Lg,\ S^\theta=g^{-1}S\}.$$

The normalizer of $\hat{G}$ in $\Aut(\Gamma)$ was given by Zhou and Feng~\cite{Zhou1}.

\begin{prop} \label{prop=Normalizer}{\rm(\cite[Theorem~1.1]{Zhou1})}
Let $\Gamma=\BiCay(G,R,L,S)$ be a connected bi-Cayley graph over a group $G$. If $I=\emptyset$ then $\N_{\rm Aut(\Gamma)}(\hat{G})=\hat{G}\rtimes F$,
and if $I\neq \emptyset$, then $\N_{\rm Aut(\Gamma)}(\hat{G})=\hat{G}\langle F,\delta_{\theta,x,y}\rangle$ for some $\delta_{\theta,x,y}\in I$.
\end{prop}

Write $N=\N_{\rm Aut(\Gamma)}(\hat{G})$ in Proposition~\ref{prop=Normalizer}.  Then $N_{1_0}=F$ and $N_{1_01_1}=\{\sigma_{\theta,1}\ |\ \theta\in \Aut(G)\ s.t.\ R^\theta=R,\ L^\theta=L,\ S^\theta=S\}$.  This implies that $F$ is faithful on $R\cup L\cup S$ because $G=\langle R,L,S\rangle$.

Recall that $G_{\a,\b,\gamma}$ is a non-abelian split metacyclic group defined by

\parbox{8cm}{
\begin{eqnarray*}
&& G_{\a,\b,\gamma}=\langle a,b\ |\ a^{p^\a}=1,\ b^{p^\b}=1,\ b^{-1}ab=a^{1+p^\gamma}\rangle,
\end{eqnarray*}}\hfill
\parbox{1cm}{
\begin{eqnarray}
\label{eq1}
\end{eqnarray}}

\f where $0<\gamma<\a\leq \b+\gamma$.

\begin{lem}\label{maintheoren-application0} Let $G_{\a,\b,\g}$ be defined as above. Then
\begin{enumerate}
\item[$1$.]$(b^{j}a^{i})^{k}=b^{kj}a^{i[(1+p^\gamma)^{(k-1)j}+(1+p^\gamma)^{(k-2)j}
+\cdots+(1+p^\gamma)^{j}+1]}$ for any $i\in \mathbb{Z}_{p^\a}$, $j\in \mathbb{Z}_{p^\b}$ and positive integer $k$;
\item[$2$.]If there is an automorphism $\theta$ of $G_{\a,\b,\g}$ such that $a^{\theta}=b^ma^n$ with $(m,p)=1$, then $\b<\a$.
\end{enumerate}
\end{lem}

\demo
With $b^{-1}ab=a^{1+p^{\gamma}}$, it is easy to prove that for any  $i\in \mathbb{Z}_{p^\a}$ and $j\in \mathbb{Z}_{p^\b}$, we have
$a^{i}b^{j}=b^{j}a^{i(1+p^\gamma)^{j}}$ and then part~1 follows by induction on $k$.

Since $\theta\in\Aut(G_{\a,\b,\gamma})$ and $a^\theta=b^ma^n$, we have  $o(b^{m}a^{n})=o(a)=p^\a$, and $\langle a\rangle \unlhd G_{\a,\b,\g}$ implies $\langle b^ma^n\rangle \unlhd G_{\a,\b,\g}$. Since $(p,m)=1$, we have $G_{\a,\b,\g}=\langle a,b^{m}a^{n}\rangle=\langle a\rangle \langle b^{m}a^{n}\rangle$, and hence $p^{\a+\b}=|G_{\a,\b,\g}|\leq p^\a \cdot p^\a$, that is, $\b\leq \a$. If $\b=\a$, then $G_{\a,\b,\g}=\langle a\rangle \times \langle b^{m}a^{n}\rangle$, contrary to the fact that $G_{\a,\b,\g}$ is non-abelian. Thus, $\b<\a$ and part~2 follows. \hfill\qed

A graph $\Gamma$ is called {\em locally transitive} if the stabilizer  $\Aut(\Gamma)_u$ for any $u\in V(\Gamma)$ is transitive on the neighborhood of $u$ in $V(\Gamma)$.

\begin{theorem}\label{maintheoren-application1}
There is no connected locally transitive bipartite bi-Cayley graphs of valency less than $p$ over a non-abelian metacyclic $p$-group $G$ for an odd prime $p$.
\end{theorem}

\demo Suppose to the contrary that $\Gamma$ is a connected locally transitive bipartite bi-Cayley graph over $G$ with valency less than $p$. Since $p$ is odd, the two orbits of $G$ are exactly the two partite sets of $\Gamma$, and we may assume that $\Gamma=\BiCay(G,\emptyset,\emptyset,S)$, where $1\in S$, $|S|<p$ and $G=\langle S\rangle$.

Let $A=\Aut(\Gamma)$. Since $\Gamma$ has valency less than $p$, $A_{1_0}$ is a $p'$-group, and by Theorem~\ref{maintheorem}, $\hat{G}\unlhd A$. Then  Proposition~\ref{prop=Normalizer} implies that $A_{1_0}=F=\{\sigma_{\theta,g}\ |\ \theta\in \Aut(G), S^\theta=g^{-1}S\}$. Note that $F$ is a group with operation $\sigma_{\theta,x}\sigma_{\delta,y}=\sigma_{\theta\delta,yx^\delta}$ for any $\sigma_{\theta,x},\sigma_{\delta,y}\in F$.
Since $1\in S$, we have $F=\{\sigma_{\theta,s}\ |\ \theta\in \Aut(G),  s\in S, S^\theta=s^{-1}S\}$, and since $\Gamma$ is locally transitive, $F$ is transitive on $S_1=\{s_1\ |\ s\in S\}$.

Set $L=\{ \theta\ |\ \sigma_{\theta,s}\in F\}$. Since $F$ is a group, $L$ is a group and the map  $\varphi: \sigma_{\theta,s}\mapsto \theta$ defines a homomorphism from $F$ to $L$. Let $K$ be the kernel of $\varphi$. If $K\not=1$ then $\sigma_{1,s}\in F$ for some $1\not=s\in S$ ($\sigma_{1,1}$ is the identity of $F$).
Note that $1_{1}^{\sigma_{1,s}}=s_1$ and $(s_{1}^{l-1})^{\sigma_{1,s}}=s_{1}^{l}$ for any positive integer $l$.
It follows that $\{s_1,s^{2}_1,s^{3}_1,\ldots\}\subseteq S_1$, and since $s$ has order $p$-power, $|S_1|\geq p$, contrary to the fact $|S_1|=|S|<p$. Thus, $K=1$ and $F\cong L$.

Assume that $G$ is non-split. By Lindenberg~\cite{Linden}, the automorphism group of $G$ is a $p$-group. Thus, $p  \mid |L|$ and hence $p\mid |A_{1_0}|$, a contradiction.

Assume that $G$ is split. Then $G=G_{\a,\b,\gamma}$, as defined in Eq~(\ref{eq1}). Since $F$ is a $p'$-group and $F\cong L$,
Proposition~\ref{prop=$p-1$-subgroup} implies that $F$ is cyclic and $|F|\mid (p-1)$. Let $|S_1|=k$ and $F=\langle \sigma_{\theta,s}\rangle$, where $\theta\in \Aut(G)$, $s\in S$ and $S^\theta=s^{-1}S$. Since $F$ is transitive on $S_1$, $\sigma_{\theta,s}$ permutes all elements in $S_1$ cyclically, and so $\sigma_{\theta,s}^{k}$ fixes all elements in $S_1$. By Proposition~\ref{prop=Normalizer}, $F$ is faithful on $S_1$, implying that $\sigma_{\theta,s}^{k}=1$. It follows that $\sigma_{\theta,s}$ has order $k$ and is regular on $S_1$. Since $F\cong L$, $\theta$ also has order $k$. Furthermore,  $S_1=1_{1}^{\langle\sigma_{\theta,s}\rangle}=\{1_1,s_1,(ss^\theta)_1,\ldots, (ss^{\theta}\cdots s^{\theta^{k-2}})_1\}$ and $1_1=(ss^{\theta}\cdots s^{\theta^{k-1}})_1$, that is, $S=\{1,s,ss^\theta,\ldots,ss^{\theta}\cdots s^{\theta^{k-2}}\}$ and

\parbox{8cm}{
\begin{eqnarray*}
&&ss^{\theta}\cdots s^{\theta^{k-1}}=1.
\end{eqnarray*}}\hfill
\parbox{1cm}{
\begin{eqnarray}
\label{eq2-1}
\end{eqnarray}}

Note that for any  $\tau\in \Aut(G)$, we have $\Gamma=\BiCay(G,\emptyset,\emptyset,S)\cong\BiCay(G,\emptyset,\emptyset,S^\tau)$, where $S^{\tau}=\{1,t,tt^{\theta^\tau},\ldots,tt^{\theta^\tau}\cdots t^{(\theta^{\tau})^{k-2}}\}$ and $tt^{\theta^\tau}\cdots t^{(\theta^{\tau})^{k-1}}=1$ with $t=s^{\tau}$. By Proposition~\ref{prop=$p-1$-subgroup}, all cyclic groups of order $k$ in $\Aut(G)$ are conjugate, and so we may assume that $\theta$ is the automorphism induced by $a\mapsto a^e$, $b\mapsto b$, where $e\in \mz_{p^\a}^{*}$ has order $k$.

Let $s=b^i a^j\in G_{\a,\b,\gamma}$. Recall that $a^{i}b^{j}=b^{j}a^{i(1+p^\gamma)^{j}}$. Since $a^\theta=a^e$ and $b^{\theta}=b$, we have $ss^{\theta}\cdots s^{\theta^{k-1}}=b^{ki}a^{\ell}$ for some $\ell\in \mz_{p^\a}$. By Eq~(\ref{eq2-1}), $b^{ki}=1$, that is, $ki=0$ in $\mz_{p^\b}$. Since $k<p$, we have $i=0$ in $\mz_{p^\b}$, and hence $G=\langle S\rangle=\langle 1,a^j,a^ja^{je},\cdots,a^ja^{je}\cdots a^{je^{k-2}} \rangle\leq \langle a\rangle$, a contradiction. This completes the proof. \hfill\qed

\begin{cor}
Let $G$ be a non-abelian metacyclic $p$-group for an odd prime $p$. Then
there exist no connected semisymmetric or arc-transitive bipartite bi-Cayley graphs over $G$ with valency less than $p$.
\end{cor}

To classify connected bipartite half-arc-transitive bi-Cayley graphs of valency less than $2p$ over non-abelian metacyclic $p$-groups, we need the following lemma.

\begin{lem}\label{lem=$e-1$ reversible and $Tx=T$}
Let $e$ be an element of order $k$ {\rm ($k\geq 2$)} in $\mathbb{Z}_{p^\a}^{*}$ with $k\di (p-1)$ and $p$ a prime. Then $e^i-1\in \mz_{p^\a}^*$ for any $1\leq i<k$, and $1+e+\cdots+e^{k-1}=0$ in $\mz_{p^\a}$.

Let $T=\{0,1,1+e,\ldots,1+e+\cdots+e^{k-2}\}=\{(e-1)^{-1}(e^i-1)\ |\ i\in \mz_k\}$. Then $T\subseteq \mz_{p^\a}$, and $Tx+y=\{tx+y\ |\ t\in T\}=T$ for $x,y\in\mz_{p^\a}$ if and only if $x=e^l$ and $y=(e-1)^{-1}(e^l-1)$ for some $l\in \mathbb{Z}_k$. In particular, $Tx=T$ if and only if $x=1$.
\end{lem}

\demo Suppose $e^i-1\not\in \mz_{p^\a}^*$ for some $1\leq i<k$. Then $p\mid (e^i-1)$, and since $e$ has order $k$, we have $e^i\not=1$. Furthermore, there exist $l\in \mz_{p^\a}^*$ ($p\nmid l$) and $1\leq s<\a$ such that $e^i=1+lp^s$. It follows that   $0=(e^i)^k-1=(1+lp^s)^k-1=klp^{s}+C_{k}^{2}(lp^s)^{2}+\cdots+C_{k}^{k-1}(lp^s)^{k-1}+(lp^s)^k$, and hence $p\mid kl$. Since $k<p$, we have $p\mid l$, a contradiction. Thus, $p\nmid (e^i-1)$, that is, $e^i-1\in \mz_{p^\a}^*$. The fact $1+e+\cdots+e^{k-1}=0$ follows from $0=e^k-1=(e-1)(1+e+\cdots+e^{k-1})$ and $e-1\in \mz_{p^\a}^*$.

Clearly, $T\subseteq \mz_{p^\a}$ and $k\in\mz_{p^\a}^*$. Then $\sum_{t\in T} t=(e-1)^{-1}\sum_{i\in \mz_k} (e^i-1)=(e-1)^{-1}[(e-1)+\cdots+(e^{k-1}-1)]=(e-1)^{-1}[-k+(1+e+\cdots+e^{k-1})]=
-k(e-1)^{-1}\in \mathbb{Z}_{p^\a}^*$.

Let $Tx+y=T$ for $x,y\in\mz_{p^\a}$. Then $\sum_{t\in T}(tx+y)=\sum_{t\in T}t$, and hence $ky=(1-x)\sum_{t\in T}t=-(1-x)k(e-1)^{-1}$.
It follows $y=(e-1)^{-1}(x-1)$ because $k\in\mz_{p^\a}^*$. Then $Tx+(e-1)^{-1}(x-1)=T$ implies $x[T(e-1)+1]=T(e-1)+1$. Since $T(e-1)+1=\{e^i\ |\ i\in \mz_k\}=\langle e
\rangle$, we have $x\langle e\rangle=\langle e\rangle$ in $\mz_{p^\a}^*$. Then $x=e^l$ for some $l\in \mz_k$, and so $y=(e-1)^{-1}(e^l-1)$.

On the other hand, let $x=e^l$ and $y=(e-1)^{-1}(e^l-1)$ with $l\in \mathbb{Z}_k$. Then $Tx+y=\{e^l(e-1)^{-1}(e^i-1)+(e-1)^{-1}(e^l-1)\ |\ i\in\mz_k\}=(e-1)^{-1}\{ e^l(e^i-1)+(e^l-1)\ |\ i\in\mz_k \}=(e-1)^{-1}\{e^{i+l}-1\ |\ i\in \mz_k\}=\{ (e-1)^{-1}(e^i-1)\ |\ i\in\mz_k \}=T$.

At last, by taking $y=0$ we have that $Tx=T$ if and only if $x=1$. \hfill\qed

Let $e$ be an element of order $k\geq 2$ in $\mathbb{Z}_{p^\a}^{*}$ with $k\di (p-1)$ and $p$ a prime. For $m\in\mz_{p^{\a-\g}}^*$ ($0<\g<\a$) and $l\in\mz_k$, let us consider solutions of the following equation in $\mz_{p^\a}$:

\parbox{8cm}{
  \begin{eqnarray*}
  && e^{l}(1+p^\gamma)^{m}=[(1+p^\gamma)^{m}-x(1-e)]^{2}.
  \end{eqnarray*}}\hfill
  \parbox{1cm}{
  \begin{eqnarray}
  \label{eq3}
  \end{eqnarray}}

Since $e^{l}(1+p^\gamma)^{m}\in\mz_{p^\a}^*$, Eq~(\ref{eq3}) has a solution if and only if $e^{l}(1+p^\gamma)^{m}$ is a square in $\mz_{p^\a}^*$. Since $\mz_{p^\a}^*\cong\mz_{p^{\a-1}(p-1)}$, squares in $\mz_{p^\a}^*$ consists of the unique subgroup of order $\frac{(p-1)}{2}p^{\a-1}$ in $\mz_{p^\a}^*$, and so Eq~(\ref{eq3}) has a solution if and only if the order of $e^{l}(1+p^\gamma)^{m}$ in $\mz_{p^\a}^*$  is a divisor of $\frac{(p-1)}{2}p^{\a-1}$. Clearly, $(1+p^\gamma)^{m}$ has order $p^{\a-\g}$, and $e^{l}$ has order $\frac{k}{(k,l)}$. Thus, Eq~(\ref{eq3}) has a solution if and only if $\frac{k}{(k,l)}\di \frac{(p-1)}{2}$. In this case, if $e^{l}(1+p^\gamma)^{m}=u^2$ for some $u\in\mz_{p^\a}^*$ then $(1-e)^{-1}[(1+p^\gamma)^m\pm u]$ are the only two solutions of Eq~(\ref{eq3}) in $\mz_{p^\a}$. It is also easy to see that Eq~(\ref{eq3}) has solutions in most cases. In fact, Eq~(\ref{eq3}) has no solution if and only if $\frac{k}{(k,l)}\nmid \frac{(p-1)}{2}$ with $k\di (p-1)$
if and only if $k$ is even, $l$ is odd and $k_2=(p-1)_2$, where $k_2$ and $(p-1)_2$ are the largest $2$-powers in $k$ and $p-1$, respectively.

\medskip
\f{\bf Construction of half-arc-transitive graphs:} Let $G_{\a,\b,\gamma}$ be the group in Eq~(\ref{eq1}). Let $m\in \mz_{p^{\a-\gamma}}^{*}$ and $k\di (p-1)$ with $k\geq 2$. Choose $0\leq l<k$ such that $\frac{k}{(k,l)}\di \frac{(p-1)}{2}$.

Take $e$ as an element of order $k$ in $\mathbb{Z}_{p^\a}^{*}$ and $n$ as a solution of Eq~(\ref{eq3}). Then $n$ is determined by $m,k,l$. Let $T=\{(e-1)^{-1}(e^i-1)\ |\ i\in \mz_k\}$ and $T'=\{(e-1)^{-1}(e^i-1)(1+p^\gamma)^{m}+e^{i}n\ |\ i\in \mz_{k}\}$. Then $T,T'\subseteq \mz_{p^\a}$. Set $U=\{a^t\ |\ t\in T\}$ and $V=\{b^ma^i\ |\ i\in T'\}$.

$$\mbox{Define:}\ \ \ \ \ \ \ \ \ \ \ \ \  \Gamma_{m,k,l}^n=\BiCay(G_{\a,\b,\gamma},\emptyset,\emptyset,U\cup V).$$

Note that $e$ is an element of order $k$ in $\mz_{p^\a}^*$ given in advance. Since Eq~(\ref{eq3}) has exactly two solutions, we also write the notation $\Gamma_{m,k,l}^n$ as $\Gamma_{m,k,l}^{\pm}$.

The following result is a classification of connected bipartite half-arc-transitive bi-Cayley graphs of valency less than $2p$ over non-abelian metacyclic $p$-groups.

\begin{theorem}\label{maintheorem-application2}
Let $G$ be a non-abelian metacyclic $p$-group for an odd prime $p$, and let $\Gamma$ be a connected bipartite bi-Cayley graph over $G$ with valency less than $2p$. Then $\Gamma$ is half-arc-transitive if and only if $\Gamma\cong \Gamma_{m,k,l}^{\pm}$ with valency $2k$ and stabilizers of $\Aut(\Gamma_{m,k,l}^{\pm})$ isomorphic to $\mz_k$.
\end{theorem}

\demo The two orbits of $G$ are exactly the two partite sets of $\Gamma$, and we may assume that $\Gamma=\BiCay(G,\emptyset,\emptyset,S)$, where $1\in S$, $|S|<2p$ and $G=\langle S\rangle$. Let $A=\Aut(\Gamma)$.

\medskip
To prove the necessity, let $\Gamma$ be half-arc-transitive. Then $A_{1_0}$ has exactly two orbits on $S_1=\{s_1\ |\ s\in S\}$, say $U_1$ and $V_1$ with $1_1\in U_1$. Let $|U|=k$. Then $S=U\cup V$, $|U|=|V|$ and $|S|=2k$. Since $k<p$, the Orbit-Stabilizer theorem implies that $A_{1_0}$ is a $p'$-group. By Theorem~\ref{maintheorem}, $\hat{G}\unlhd A$, and by Proposition~\ref{prop=Normalizer}, $A_{1_0}=F=\{\sigma_{\theta,g}\ |\ \theta\in \Aut(G), S^\theta=g^{-1}S\}$.
Since $1\in S$, we have $F=\{\sigma_{\theta,s}\ |\ \theta\in \Aut(G), s\in S, S^\theta=s^{-1}S\}$. Note that $F$ is a group with operation $\sigma_{\theta,x}\sigma_{\delta,y}=\sigma_{\theta\delta,yx^\delta}$ for any $\sigma_{\theta,x}$, $\sigma_{\delta,y}\in F$.

Set $L:=\{ \theta\ |\ \sigma_{\theta,s}\in F\}$. Since $F$ is a group, $L$ is a group and the map $\varphi: \sigma_{\theta,s}\mapsto \theta$ defines a homomorphism from $F$ to $L$. Let $K$ be the kernel of $\varphi$. Suppose $K\neq 1$. Then there is $1\not=\sigma_{1,s}\in K$ with $s\not=1$, and so $1_{1}^{\langle\sigma_{1,s}\rangle}
=\{1_1,s_1,s^2_1,\ldots\}\subseteq U_1$. Since $s$ has order $p$-power, $|U_1|\geq p$ and hence $k\geq p$, a contradiction. Thus, $K=1$ and  $F\cong L$.

Assume that $G$ is non-split. By Lindenberg~\cite{Linden}, the automorphism group of $G$ is a $p$-group. Thus, $p\di|L|$ and hence $p\di |A_{1_0}|$, a contradiction.

Assume that $G$ is split. Then $G=G_{\a,\b,\gamma}$, as defined in Eq~(\ref{eq1}). Since $F$ is a $p'$-group and $F\cong L$, Proposition~\ref{prop=$p-1$-subgroup} implies that $F$ is cyclic and $|F|\di (p-1)$. Let $F=\langle \sigma_{\theta,s}\rangle$, where $\theta\in \Aut(G)$, $s\in S$ and $S^\theta=s^{-1}S$. The transitivity of $F$ implies that $\sigma_{\theta,s}$ permutes all elements in both $U_1$ and $V_1$ cyclically. Thus, $\sigma_{\theta,s}^{k}$ fixes all elements in $S_1$, and by Proposition~\ref{prop=Normalizer}, $\sigma_{\theta,s}^{k}=1$ and hence $\sigma_{\theta,s}$ has order $k$ and is regular on both $U_1$ and $V_1$. It follows that   $A_{1_0}=F=\langle \sigma_{\theta,s}\rangle\cong\mz_k$. Furthermore,  $U_1=1_{1}^{\langle\sigma_{\theta,s}\rangle}=\{1_1, s_1, (ss^\theta)_1, \ldots, (ss^{\theta}\cdots s^{\theta^{k-2}})_1\}$,
and $V_1=t_{1}^{\langle\sigma_{\theta,s}\rangle}=\{t_1, (st^\theta)_1, (ss^{\theta}t^{\theta^2})_1, \ldots, (ss^{\theta}\cdots s^{\theta^{k-2}}t^{\theta^{k-1}})_1\}$ with $(ss^{\theta}\cdots s^{\theta^{k-1}})_1=1_1$ for any $t\in V$, that is,
$U=\{1, s, ss^\theta, \ldots, ss^{\theta}\cdots s^{\theta^{k-2}}\}$, $V=\{t, st^\theta, ss^{\theta}t^{\theta^2}, \ldots, ss^{\theta}\cdots s^{\theta^{k-2}}t^{\theta^{k-1}}\}$,
$\theta$ has order $k$ with $k\di (p-1)$, and

\parbox{8cm}{
\begin{eqnarray*}
&&ss^{\theta}\cdots s^{\theta^{k-1}}=1.
\end{eqnarray*}}\hfill
\parbox{1cm}{
\begin{eqnarray}
\label{eq6-1}
\end{eqnarray}}

\f By Proposition~\ref{prop=$p-1$-subgroup}, we may assume that $\theta$ is the automorphism induced by $a\mapsto a^e$, $b\mapsto b$, where $e\in \mathbb{Z}_{p^\a}^{*}$ has order $k$.

Let $s=b^{i}a^j$ and $t=b^{m}a^{n}$ with $i,m\in \mz_{p^\b}$ and $j,n\in\mz_{p^\a}$. Then $ss^{\theta}\cdots s^{\theta^{k-1}}=b^{ki}a^{\ell}$ for some $\ell\in \mathbb{Z}_{p^\a}$. By Eq~(\ref{eq6-1}), $b^{ki}=1$, that is, $ki=0$ in $\mathbb{Z}_{p^\b}$. Since $k<p$, we have $i=0$ in $\mz_{p^\b}$, and hence $s=a^j$. Since $ss^{\theta}\cdots s^{\theta^{i-1}}=a^ja^{je}\cdots a^{je^{i-1}}=a^{j(e-1)^{-1}(e^i-1)}$, we have
$$U=\{1, a^j, a^{j}a^{je}, \ldots, a^{j}a^{je}\cdots a^{je^{k-2}}\}=\{a^{j(e-1)^{-1}(e^i-1)}\ |\ i\in \mz_k\}.$$

Since $e\in \mathbb{Z}_{p^\a}^{*}$, any element of order $k$ in $\mathbb{Z}_{p^\a}^{*}$ can be written as $e^q$ with $(q,k)=1$. By Lemma~\ref{lem=$e-1$ reversible and $Tx=T$}, $e-1\in\mz_{p^\a}^*$ and $e^q-1\in\mz_{p^\a}^*$, and so $G$ has an automorphism $\rho$ induced by $a\mapsto a^{(e-1)(e^q-1)^{-1}}$ and $b\mapsto b$. It follows that $U^\rho=\{a^{(e-1)(e^q-1)^{-1}j(e-1)^{-1}(e^i-1)}\ |\ i\in \mz_k\}=\{a^{j(e^{q}-1)^{-1}((e^q)^i-1)}\ |\ i\in \mz_k\}$. It is easy to check that $\theta^\rho=\theta$. Thus, we may take $e$ as an special element of order $k$ in $\mathbb{Z}_{p^\a}^{*}$ in advance.

Since $a^{j(e-1)^{-1}(e^i-1)}b^m=b^ma^{j(e-1)^{-1}(e^i-1)(1+p^\gamma)^{m}}$ and $(b^ma^n)^{\theta^i}=b^ma^{e^in}$, we have $$V=\{a^{j(e-1)^{-1}(e^i-1)}(b^ma^n)^{\theta^i}\ |\ i\in \mz_k\}=\{b^{m}a^{j(e-1)^{-1}(e^i-1)(1+p^\gamma)^{m}+e^{i}n} \ |\ i\in \mz_k\}.$$
Thus, $G=\langle a^j, a^n, b^m\rangle$ because $G=\langle S\rangle=\langle U\cup V\rangle\leq \langle a^j, a^n, b^m\rangle$.

Since $\Gamma$ is half-arc-transitive, Proposition~\ref{prop=Normalizer} implies that there exists  $\delta_{\lambda,x,y}\in I$ such that $(1_0,1_1)^{\delta_{\lambda,x,y}}$ $=$ $((b^{m}a^{n})_1,1_0)$ with $\lambda\in \Aut(G)$ and $S^\lambda=y^{-1}S^{-1}x$. In particular, $(b^{m}a^{n})_1=1_{0}^{\delta_{\lambda,x,y}}=x_1$ and $1_0=1_{1}^{\delta_{\lambda,x,y}}=y_0$.
It follows that $x=b^{m}a^{n}$, $y=1$ and $S^\lambda=S^{-1}b^{m}a^{n}=U^{-1}b^{m}a^{n}\cup V^{-1}b^{m}a^{n}$. Furthermore,
$$U^{-1}b^{m}a^{n}=\{a^{-j(e-1)^{-1}(e^i-1)}b^{m}a^{n}\ |\ i\in \mz_k\}=\{b^{m}a^{-j(e-1)^{-1}(e^i-1)(1+p^\gamma)^{m}+n}\ |\ i\in \mz_k\},$$
and since $(b^{m}a^{j(e-1)^{-1}(e^i-1)(1+p^\gamma)^{m}+e^{i}n})^{-1}b^{m}a^{n}=a^{-j(e-1)^{-1}(e^i-1)(1+p^\gamma)^{m}+n(1-e^i)}$, we have
$$V^{-1}b^{m}a^{n}=\{a^{-j(e-1)^{-1}(e^i-1)(1+p^\gamma)^{m}+n(1-e^i)}\ |\ i\in \mz_k\}.$$

Suppose $p\di j$. Since $G=\langle a^j, a^n, b^m\rangle$, we have $p\nmid n$ and $p\nmid m$. By Proposition~\ref{prop=element order}, every element in both $V$ and $U^{-1}b^{m}a^{n}$ has order $\max\{p^\a,p^\b\}$. Clearly, every element in $U$ has order less than $p^\a$, but the element $a^{-j(1+p^\gamma)^m+n(1-e)}\in V^{-1}b^{m}a^{n}$ has order $p^\a$ because $p\nmid (1-e)$ by Lemma~\ref{lem=$e-1$ reversible and $Tx=T$}. This is impossible as $\lambda\in \Aut(G)$ and $(U\cup V)^\lambda=S^\lambda=S^{-1}b^{m}a^{n}=U^{-1}b^{m}a^{n}\cup V^{-1}b^{m}a^{n}$. Thus, $p\nmid j$. Furthermore, $p\nmid m$ and so $m\in \mz_{p^\b}^*$.

Now, there is an automorphism of $G$ mapping $a^j$ to $a$ and $b$ to $b$, and so we may assume $j=1$ and $s=a$. It follows that

\parbox{8cm}{
\begin{eqnarray*}
&& U=\{a^\eta\ |\ \eta\in T\}, \ \ \ T=\{(e-1)^{-1}(e^i-1)\ |\ i\in \mz_k\};\\
&& V=\{b^ma^\eta\ |\ \eta\in T'\}, \ \ \ T'=\{(e-1)^{-1}(e^i-1)(1+p^\gamma)^{m}+e^{i}n\ |\ i\in \mz_{k}\}.
\end{eqnarray*}}\hfill
\parbox{1cm}{
\begin{eqnarray}
\label{eq-T}\\ \label{eq-T'}
\end{eqnarray}}

\f As $(e-1)^{-1}(e^i-1)(1+p^\gamma)^{m}+e^{i}n=[(e-1)^{-1}(e^i-1)][(1+p^\gamma)^m+n(e-1)]+n$, we have

\parbox{8cm}{
\begin{eqnarray*}
&& T'=T [(1+p^\gamma)^m+n(e-1)]+n.
\end{eqnarray*}}\hfill
\parbox{1cm}{
\begin{eqnarray}
\label{eqT'-T}
\end{eqnarray}}

\f Since $-(e-1)^{-1}(e^i-1)(1+p^\gamma)^{m}+n=[(e-1)^{-1}(e^i-1)][-(1+p^\gamma)^{m}]+n$ and $-(e-1)^{-1}(e^i-1)(1+p^\gamma)^{m}+n(1-e^{i})=[(e-1)^{-1}(e^i-1)][-(1+p^\gamma)^{m}+n(1-e)]$,
we have

\parbox{8cm}{
\begin{eqnarray*}
&& U^{-1}b^{m}a^{n}=\{b^{m}a^\eta\ |\ \eta\in T_1\},\ \ \ T_1=T [-(1+p^\gamma)^{m}]+n;\\
&& V^{-1}b^{m}a^{n}=\{a^\eta\ |\ \eta\in T_1'\},\ \ \ T_1'=T[-(1+p^\gamma)^{m}+n(1-e)].
\end{eqnarray*}}\hfill
\parbox{1cm}{
\begin{eqnarray}
\label{eq-T1}\\ \label{eq-T1'}
\end{eqnarray}}

Noting that $T,T',T_1,T_1'\subseteq \mz_{p^\a}$, we have $U,V,U^{-1}b^{m}a^{n},V^{-1}b^{m}a^{n}\subseteq G$.

\medskip
\f{\bf Claim:} $a^\lambda\in V^{-1}b^{m}a^{n}$.

Suppose $a^\lambda\notin V^{-1}b^{m}a^{n}$. Since $a^\lambda\in S^{\lambda}=U^{-1}b^{m}a^{n}\cup V^{-1}b^{m}a^{n}$, we have $a^\lambda\in U^{-1}b^{m}a^{n}$, that is, $a^\lambda=b^{m}a^{\mu}$ for $\mu\in T_{1}$, implying that $\b<\a$ by Lemma~\ref{maintheoren-application0}. Clearly, $k\geq 2$ as $\Gamma$ is half-arc-transitive.

Let $k>2$. Then $(a^{1+e})^\lambda=(b^{m}a^{\mu})^{1+e}\in U^{-1}b^{m}a^{n}$. Considering the power of $b$, we have $m(1+e)=m$ in $\mz_{p^\b}$ as $p\nmid m$ and $p\nmid (1+e)$ by Lemma~\ref{lem=$e-1$ reversible and $Tx=T$}. It follows that $e=0$ in $\mz_{p^\b}$, which is impossible because $e\in \mz_{p^\a}^{*}$ implies $(p,e)=1$.

Let $k=2$. Then $T=\{0,1\}$ and $e=-1$ in $\mz_{p^\a}$. By Eqs~(\ref{eq-T}) and (\ref{eq-T'}), $S$ $=$ $\{1$, $a$, $b^{m}a^{n}$, $b^{m}a^{(1+p^\gamma)^{m}-n}\}$, and by Eqs~(\ref{eq-T1}) and (\ref{eq-T1'}), $S^{-1}b^{m}a^{n}$ $=$ $\{1$, $a^{-(1+p^\gamma)^{m}+2n}$, $b^{m}a^{n},b^{m}a^{-(1+p^\gamma)^{m}+n}\}$. Note that $a^\lambda\in U^{-1}b^{m}a^{n}=\{b^{m}a^{n},b^{m}a^{-(1+p^\gamma)^{m}+n}\}$.

\medskip
\f{\bf Case 1:}  $a^\lambda=b^{m}a^{n}$.

As $S^{\lambda}=S^{-1}b^{m}a^{n}$, we have $((b^{m}a^{n})^\lambda, (b^{m}a^{(1+p^\gamma)^{m}-n})^\lambda)=(a^{-(1+p^\gamma)^{m}+2n},b^{m}a^{-(1+p^\gamma)^{m}+n})$ or $(b^{m}a^{-(1+p^\gamma)^{m}+n},a^{-(1+p^\gamma)^{m}+2n})$. For the former,
$b^{m}a^{-(1+p^\gamma)^{m}+n}=(b^{m}a^{(1+p^\gamma)^{m}-n})^\lambda$ $=$  $[(b^{m}a^{n})a^{(1+p^\gamma)^{m}-2n}]^\lambda$ $=$ $a^{-(1+p^\gamma)^{m}+2n}$ $(b^{m}a^{n})^{(1+p^\gamma)^{m}-2n}$, implying that $m=m[(1+p^\gamma)^{m}-2n]$ in $\mz_{p^\b}$, and since $p\nmid m$, we have $p\di n$. This is impossible because otherwise $p^\a=o(a^\lambda)=o(b^{m}a^{n})<p^\a$ ($\b<\a$). For the latter,
$a^{-(1+p^\gamma)^{m}+2n}$ $=$ $(b^{m}a^{(1+p^\gamma)^{m}-n})^\lambda$ $=$ $[(b^{m}a^{n})a^{(1+p^\gamma)^{m}-2n}]^\lambda$ $=$ $b^{m}a^{-(1+p^\gamma)^{m}+n}(b^{m}a^{n})^{(1+p^\gamma)^{m}-2n}$. Thus,  $0=m+m[(1+p^\gamma)^{m}-2n]$ in $\mz_{p^\b}$, and hence $p\di (1-n)$, but it is also impossible because otherwise $p^\a=o(a^\lambda)=o(b^{m}a^{n})=o((b^{m}a^{n})^\lambda)=o(b^{m}a^{-(1+p^\gamma)^{m}+n})<p^\a$.

\medskip
\f{\bf Case 2:} $a^\lambda=b^{m}a^{-(1+p^\gamma)^{m}+n}$.

As $S^{\lambda}=S^{-1}b^{m}a^{n}$, we have $((b^{m}a^{n})^\lambda, (b^{m}a^{(1+p^\gamma)^{m}-n})^\lambda)=(a^{-(1+p^\gamma)^{m}+2n},b^{m}a^{n})$ or $(b^{m}a^{n}$, $a^{-(1+p^\gamma)^{m}+2n})$. For the former,
$b^{m}a^{n}=(b^{m}a^{(1+p^\gamma)^{m}-n})^\lambda$ $=$  $[(b^{m}a^{n})a^{(1+p^\gamma)^{m}-2n}]^\lambda$ $=$ $a^{-(1+p^\gamma)^{m}+2n}$ $(b^{m}a^{-(1+p^\gamma)^{m}+n})^{(1+p^\gamma)^{m}-2n}$, implying that $m=m[(1+p^\gamma)^{m}-2n]$ in $\mz_{p^\b}$, and since $p\nmid m$, we have $p\di n$. This is impossible because otherwise $p^\a=o(a^\lambda)=o(b^{m}a^{-(1+p^\gamma)^{m}+n})=o(b^{m}a^{(1+p^\gamma)^{m}-n})
=o((b^{m}a^{(1+p^\gamma)^{m}-n})^\lambda)=o(b^{m}a^{n})<p^\a$ ($\b<\a$).
For the latter,
$a^{-(1+p^\gamma)^{m}+2n}$ $=$ $(b^{m}a^{(1+p^\gamma)^{m}-n})^\lambda$ $=$ $[(b^{m}a^{n})a^{(1+p^\gamma)^{m}-2n}]^\lambda$ $=$ $b^{m}a^{n}(b^{m}a^{-(1+p^\gamma)^{m}+n})^{(1+p^\gamma)^{m}-2n}$. Thus,  $0=m+m[(1+p^\gamma)^{m}-2n]$ in $\mz_{p^\b}$, and hence $p\di (1-n)$, but it is also impossible because otherwise $p^\a=o(a^\lambda)=o(b^{m}a^{-(1+p^\gamma)^{m}+n})<p^\a$.
This completes the proof of Claim.
\medskip

By Claim, $a^\lambda=a^\mu\in V^{-1}b^{m}a^{n}$ for some $\mu\in T_{1}^{'}$. Since $p^\a=o(a^\lambda)=o(a^\mu)$, we have $\mu\in \mz_{p^\a}^{*}$. By Eqs~(\ref{eq-T}) and (\ref{eq-T1'}), $U^\lambda=\{a^{\eta\mu}\ |\ \eta\in T\}=\{a^\eta\ |\ \eta\in T\mu\}\subseteq V^{-1}b^{m}a^{n}$. Then $U^\lambda=V^{-1}b^{m}a^{n}=\{a^\eta\ |\ \eta\in T_{1}^{'}\}$, and so $T\mu=T_{1}^{'}$ in $\mz_{p^\a}$. By Eq~(\ref{eq-T1'}), $T\mu=T[-(1+p^\gamma)^{m}+n(1-e)]$. Since $p\nmid \mu$, we have $T=T[-(1+p^\gamma)^{m}+n(1-e)]\mu^{-1}$. By Lemma~\ref{lem=$e-1$ reversible and $Tx=T$}, $\mu=-(1+p^\gamma)^{m}+n(1-e)$.

Since $S^\lambda=S^{-1}b^{m}a^{n}=U^{-1}b^{m}a^{n}\cup V^{-1}b^{m}a^{n}$, we have $V^\lambda=U^{-1}b^{m}a^{n}$. In particular, $(b^{m}a^{n})^\lambda=b^{m}a^{\nu}$ for some $\nu\in T_{1}$.
For $\eta\in T'$, since $(b^{m}a^{\eta})^\lambda=[(b^{m}a^{n})a^{\eta-n}]^{\lambda}=b^{m}a^{\nu}a^{\eta\mu- \mu n}=b^{m}a^{\eta\mu-\mu n +\nu}$, we have that
$\{b^{m}a^{\eta}\ |\ \eta\in T_1\}=U^{-1}b^{m}a^{n}=V^\lambda=\{(b^{m}a^{\eta})^\lambda\ |\ \eta\in T'\}=\{b^{m}a^{\eta\mu-\mu n +\nu}\ |\ \eta\in T'\}=\{b^{m}a^{\eta}\ |\ \eta\in T^{'}\mu-\mu n+\nu\}$. By Eqs~(\ref{eqT'-T}) and (\ref{eq-T1}),
$T[-(1+p^\gamma)^{m}]+n=T_1=T'\mu-\mu n+\nu=T[(1+p^\gamma)^{m}+n(e-1)]\mu+\mu n-\mu n+\nu$ in $\mz_{p^\a}$.
Thus,    $T[(1+p^\gamma)^{m}-n(1-e)]^{2}(1+p^\gamma)^{-m}-(\nu-n)(1+p^\gamma)^{-m}=T$.
By Lemma~\ref{lem=$e-1$ reversible and $Tx=T$}, there exists $l\in\mathbb{Z}_k$ such that $e^l=[(1+p^\gamma)^{m}-n(1-e)]^{2}(1+p^\gamma)^{-m}$,
that is, $n$ satisfies Eq~(\ref{eq3}).

Recall that  $\a-\gamma\leq \b$ and $m\in \mz_{p^\b}^{*}$. Let $m=m_1+lp^{\a-\gamma}$ with $m_1\in \mz_{p^{\a-\gamma}}^{*}$. Since $(1+p^\gamma)$ has order $p^{\a-\gamma}$ in $\mz_{p^\a}$, we have $(1+p^\gamma)^{m}=(1+p^\gamma)^{m_1+lp^{\a-\gamma}}=(1+p^\gamma)^{m_1}$ in $\mz_{p^\a}$. This implies that replacing  $m$ by $m_1$, Eq~(\ref{eq3}) has the same solutions, and $T'=\{(e-1)^{-1}(e^i-1)(1+p^\gamma)^{m_1}+e^{i}n\ |\ i\in \mz_{k}\}\subseteq \mz_{p^\a}$. The automorphism of $G$ induced by  $a\mapsto a$ and $b\mapsto b^{m_{1}m^{-1}}$, maps $U$ to $U$, and $V=\{b^ma^\eta\ |\ \eta\in T'\}$ to $\{b^{m_1}a^\eta\ |\ \eta\in T'\}$.  Thus, we may assume that $m\in \mz_{p^{\a-\gamma}}^{*}$, and therefore, $\Gamma\cong\Gamma_{m,k,l}^n$.

\medskip
 We now prove the sufficiency.  Let $p$ be an odd prime and $k$ a positive integer with $k\di (p-1)$ and $k\geq 2$. Let $l\in \mz_k$, $m\in \mz_{p^{\a-\gamma}}^{*}$ and $n$ satisfies Eq~(\ref{eq3}). Take $e$ as an element of order $k$ in $\mz_{p^\a}^{*}$.  Let $S=U\cup V$, where $U=\{a^{\eta}\ |\ \eta\in T\}$ with $T=\{(e-1)^{-1}(e^i-1)\ |\ i\in \mz_k\}$ and $V=\{b^{m}a^{\eta}\ |\ \eta\in T'\}$ with $T'=\{(e-1)^{-1}(e^i-1)(1+p^\gamma)^{m}+e^{i}n\ |\ i\in \mz_k\}$.
Denote by $\Gamma_{m,k,l}^{n}$ the bi-Cayley graph $\BiCay(G_{\a,\b,\g},\emptyset,\emptyset,S)$, as constructed before Theorem~\ref{maintheorem-application2}. Clearly, $1\in U$ and $G_{\a,\b,\g}=\langle S\rangle$, implying that $\Gamma_{m,k,l}^{n}$ is connected.
Note that $$T'=T [(1+p^\gamma)^m+n(e-1)]+n.$$
To finish the proof, we only need to show that $\Gamma_{m,k,l}^{n}$ is half-arc-transitive.  Let $A=\Aut(\Gamma_{m,k,l}^{n})$ and $G=G_{\a,\b,\gamma}$.

We first claim $p\nmid |A_{1_0}|$. Suppose to the contrary that $p\mid |A_{1_0}|$. Let $P$ is a Sylow $p$-subgroup of $A$ containing  $\hat{G}$ and let $X=\N_{A}(\hat{G})$. Then $\hat{G}<P$, and hence $\hat{G}< \N_{P}(\hat{G})\leq X$. In particular, $p\di |X:\hat{G}|$, and so $p\di |X_{1_0}|$. By Proposition~\ref{prop=Normalizer}, $X_{1_0}=F=\{\sigma_{\theta,g}\ |\ \theta\in \Aut(G), S^\theta=g^{-1}S\}=\{\sigma_{\theta,s}\ |\ \theta\in \Aut(G), s\in S, S^\theta=s^{-1}S\}$ as $1\in S$.

Let $\tau$ be the automorphism of $G$ induced by $a\mapsto a^{e}$ and $b\mapsto b$.
We now prove $\sigma_{\tau,a}\in X_{1_0}$, which is equivalent to show that $S^\tau=a^{-1}S$. Clearly, $U^\tau=\{a^{e\eta}\ |\ \eta\in T\}=\{a^{\eta}\ |\ \eta\in Te\}$ and $a^{-1}U=\{a^{\eta-1}\ |\ \eta\in T\}=\{a^{\eta}\ |\ \eta\in T-1\}$. By Lemma~\ref{lem=$e-1$ reversible and $Tx=T$}, $Te=T-1$ and hence  $U^\tau= a^{-1}U$. Similarly,  $V^\tau=\{b^{m}a^{e\eta}\ |\ \eta\in T'\}=\{b^{m}a^{\eta}\ |\ \eta\in T'e\}$ and $a^{-1}V=\{a^{-1}b^{m}a^{\eta}\ |\ \eta\in T'\}=\{b^ma^{-(1+p^\g)^m}a^{\eta}\ |\ \eta\in T'\}=\{b^ma^{\eta}\ |\ \eta\in T'-(1+p^\g)^m\}$. By Eq~(\ref{eq3}),  $(1+p^\gamma)^m+n(e-1)\in\mz_{p^\a}^*$, and hence $Te=T-1$ implies $T [(1+p^\gamma)^m+n(e-1)]e+ne=T [(1+p^\gamma)^m+n(e-1)]+n-(1+p^\g)^m$. Since $T'=T [(1+p^\gamma)^m+n(e-1)]+n$, we have $T'e=T'-(1+p^\g)^m$, that is, $V^\tau=a^{-1}V$. It follows that $S^\tau=a^{-1}S$, as required.

Note that $U_1=1_{1}^{\langle\sigma_{\tau,a}\rangle}$ and $V_1=(b^{m}a^{n})_{1}^{\langle\sigma_{\tau,a}\rangle}$. Then  $X_{1_0}$ either has two orbits of length $k$ on $S_1$, or is transitive on $S_1$.
Since $p\di |X_{1_0}|$ and $X_{1_0}$ acts faithfully on $S_1$ by Proposition~\ref{prop=Normalizer}, any element of order $p$  of  $X_{1_0}$ has an orbit of length $p$ on $S_1$, and since $k<p$,  $X_{1_0}$ is transitive on $S_1$.  From $|X_{1_0}|=|X_{1_{0} 1_{1}}||1_{1}^{X_{1_0}}|=|X_{1_{0} 1_{1}}|\cdot 2k$, we have $p\mid |X_{1_{0}1_{1}}|$. By Proposition~\ref{prop=Normalizer}, $X_{1_{0}1_{1}}=\{\sigma_{\theta,1}\ |\ \theta\in \Aut(G), S^\theta=S\}$. Let $\sigma_{\theta,1}\in X_{1_{0}1_{1}}$ be of order $p$ with $\theta\in \Aut(G)$. Then $\theta$ has order $p$ and $S^\theta=S$. Recall that $k\geq 2$.

Assume $k>2$. Since $S^\theta=S$, we have $a^\theta\in S=U\cup V$. If $a^\theta\in V$ then $a^\theta=b^{m}a^{i}$ for some $i\in T'$. Since $k>2$,  $(a^{1+e})^\theta=(b^{m}a^{i})^{1+e}\in V$, and by considering the powers of $b$, we have $m(1+e)=m$ in $\mz_{p^\b}$ because $p\nmid (1+e)$ by Lemma~\ref{lem=$e-1$ reversible and $Tx=T$}. Since $m\in\mz_{p^{\a-\g}}^*$, we have $p\nmid m$ and so $e=0$ in $\mz_{p^\b}$. It follows that $p\di e$, contrary to the fact that $e\in \mz_{p^\a}^{*}$. Thus, $a^\theta\in U$, and hence, $a^\theta=a^j$ for some $j\in T$. If $a^\theta\neq a$ then $a_{1}^{\sigma_{\theta,1}}=\{a_1,a_{1}^\theta,\ldots, a_{1}^{\theta^{p-1}}\}$ is an orbit of length $p$ of $\sigma_{\theta,1}$ on $S_1$, which is impossible because there are exactly $k<p$ elements of type $a^j$ in $S$.
Thus, $a^{\theta}=a$ and $\theta$ fixes $U$ pointwise. Furthermore, $\theta$ also fixes $V$ pointwise because $|V|=k<p$. It follows that $\theta=1$ as $G=\langle S\rangle$, a contradiction.

Assume $k=2$. Then $e=-1$ in $\mz_{p^\a}$ and
$S_1=\{1_1,a_1,(b^{m}a^{n})_1,(b^{m}a^{(1+p^\gamma)^{m}-n})_1\}$. Since $1_{1}^{\sigma_{\theta,1}}=1_1$ and $p\geq 3$, $\sigma_{\theta,1}$ has order $3$ and we may assume that $a_{1}^{\sigma_{\theta,1}}=(b^{m}a^{n})_1$, $(b^{m}a^{n})_{1}^{\sigma_{\theta,1}}=(b^{m}a^{(1+p^\gamma)^{m}-n})_1$ and $(b^{m}a^{(1+p^\gamma)^{m}-n})_{1}^{\sigma_{\theta,1}}=a_1$  (replace  $\sigma_{\theta,1}$ by $\sigma_{\theta,1}^2$ if necessary), that is, $a^\theta=b^{m}a^{n}$, $(b^{m}a^{n})^{\theta}=b^{m}a^{(1+p^\gamma)^{m}-n}$ and $(b^{m}a^{(1+p^\gamma)^{m}-n})^\theta=a$. By Corollary~\ref{maintheoren-application0}, $\b<\a$.
It follows that   $a=(b^{m}a^{(1+p^\gamma)^{m}-n})^\theta=[(b^{m}a^{n})a^{(1+p^\gamma)^{m}-2n}]^\theta
=b^{m}a^{(1+p^\gamma)^{m}-n}(b^{m}a^{n})^{(1+p^\gamma)^{m}-2n}$, and so  $0=m+m[(1+p^\gamma)^{m}-2n]$ in $\mz_{p^\b}$. Thus, $p\di (1-n)$, which is impossible because otherwise $p^\a=o(a^{\theta^2})=o(b^{m}a^{(1+p^\gamma)^{m}-n})<p^\a$. Thus, $p\nmid |A_{1_0}|$, as claimed.

By Theorem~\ref{maintheorem}, $\Gamma_{m,k,l}^{n}$ is normal, and by Proposition~\ref{prop=Normalizer}, $A_{1_0}=X_{1_0}=F=\{\sigma_{\theta,s}\ |\ \sigma\in \Aut(G),s\in S, S^\theta=s^{-1}S\}$. Set $L=\{\theta\ |\ \sigma_{\theta,s}\in F\}$. Then $F\cong L\leq \Aut(G)$, and Proposition~\ref{prop=$p-1$-subgroup} implies that $F\lesssim \mz_{p-1}$. Noting that
$\sigma_{\tau,a}\in A_{1_0}$, $F$ is transitive on $S_1$ or has two orbits. Since $F$ is cyclic and faithful on $S_1$, $F$ is regular on $S_1$ for the former, and $F=\langle \sigma_{\tau,a}\rangle$ for the latter.

Suppose that $F$ is regular on $S_1$. Then $F$ is an cyclic group of order $2k$, and $|F:\langle \sigma_{\tau,a}\rangle|=2$. Thus, $\langle \sigma_{\tau,a}\rangle\unlhd F$, and so $F$ interchanges the two orbits $U_1$ and $V_1$ of $\langle \sigma_{\tau,a}\rangle$. By the regularity of $F$, there exists a $\sigma_{\theta,s}\in F$ such that $1_{1}^{\sigma_{\theta,s}}=(b^{m}a^{n})_1$, which implies that $s=b^{m}a^{n}$ and $S^\theta=s^{-1}S=(b^{m}a^{n})^{-1}S$. Since $F$ interchanges $U_1$ and $V_1$, we have $a^\theta\in (b^{m}a^{n})^{-1}V$.
It is easy to see that $(b^{m}a^{n})^{-1}S=(b^{m}a^{n})^{-1}U\cup (b^{m}a^{n})^{-1}V$, where $(b^{m}a^{n})^{-1}U=\{(b^{m}a^{n})^{-1}a^\eta\ |\ \eta\in T\}=\{b^{-m}a^{-n(1+p^\gamma)^{-m}+\eta}\ |\ \eta\in T\}=\{b^{-m}a^{\eta}\ |\ \eta\in T-n(1+p^\gamma)^{-m}\}$ and $(b^{m}a^{n})^{-1}V=\{(b^{m}a^{n})^{-1}b^{m}a^{\eta}\ |\ \eta\in T'\}=\{a^{-n+\eta}\ |\ \eta\in T'\}=\{a^{\eta}\ |\ \eta\in T'-n\}$. Since $T'=T[(1+p^\gamma)^{m}+n(e-1)]+n$, we have $(b^{m}a^{n})^{-1}V=\{a^{\eta}\ |\ \eta\in T[(1+p^\gamma)^{m}+n(e-1)]\}$.

Let $a^\theta=a^r\in (b^{m}a^{n})^{-1}V$ for some $r\in T[(1+p^\gamma)^{m}+n(e-1)]$. Since $p^\a=o(a^\theta)=o(a^r)$, we have $r\in \mz_{p^\a}^{*}$. Note that $U^\theta=\{a^{\eta r}\ |\ \eta\in T\}=\{a^{\eta}\ |\ \eta\in Tr\}\subseteq (b^{m}a^{n})^{-1}V$. Then $U^\theta=(b^{m}a^{n})^{-1}V=\{a^{\eta}\ |\ \eta\in T[(1+p^\gamma)^{m}+n(e-1)]\}$, and so $Tr=T[(1+p^\gamma)^{m}+n(e-1)]$ in $\mz_{p^\a}$. By Lemma~\ref{lem=$e-1$ reversible and $Tx=T$}, $r=(1+p^\gamma)^{m}+n(e-1)$.

Since $S^\theta=(b^{m}a^{n})^{-1}S=(b^{m}a^{n})^{-1}U\cup (b^{m}a^{n})^{-1}V$, we have $V^\theta=(b^{m}a^{n})^{-1}U$. In particular, $(b^{m}a^{n})^\theta=b^{-m}a^{t}$ for some $t\in T-n(1+p^\gamma)^{-m}$. For $\eta\in T'$, since $(b^{m}a^{\eta})^\theta=[(b^{m}a^{n})a^{\eta-n}]^\theta
=b^{-m}a^{t}a^{r(\eta-n)}=b^{-m}a^{r\eta-rn+t}$, we have $\{b^{-m}a^{\eta}\ |\ \eta\in T-n(1+p^\gamma)^{-m}\}=(b^{m}a^{n})^{-1}U=V^\theta=\{(b^{m}a^{\eta})^\theta\ |\ \eta\in T'\}=\{b^{-m}a^{r\eta-rn+t}\ |\ \eta\in T'\}=\{b^{-m}a^\eta\ |\ \eta\in T^{'}r-rn+t\}$. This implies that $T-n(1+p^\gamma)^{-m}=T^{'}r-rn+t=Tr[(1+p^\g)^{m}+n(e-1)]+rn-rn+t
=Tr[(1+p^\g)^{m}+n(e-1)]+t=T[(1+p^\g)^{m}+n(e-1)]^{2}+t$ in $\mz_{p^\a}$.
By definition of $\Gamma_{m,k,l}^{n}$,  Eq~(\ref{eq3}) holds, that is, $e^{l}(1+p^\g)^{m}=[(1+p^\g)^{m}+n(e-1)]^{2}$. It follows that  $T-n(1+p^\gamma)^{-m}=Te^{l}(1+p^\g)^{m}+t$, and hence $T=Te^{l}(1+p^\g)^{m}+t+n(1+p^\gamma)^{-m}$. By Lemma~\ref{lem=$e-1$ reversible and $Tx=T$}, there exists $l'\in \mz_k$ such that
$e^{l'}=e^{l}(1+p^\g)^{m}$ in $\mz_{p^\a}$, that is, $e^{l'-l}=(1+p^\g)^{m}$. Since $e$ is an element of order $k$, we have $(1+p^\g)^{mk}=1$ in $\mz_{p^\a}$ and since $(mk,p)=1$, we have $p^\g=0$  in $\mz_{p^\a}$, implying that $\a\di\g$, which is impossible because $0<\g<\a$.

Thus, $F$ cannot be regular on $S_1$, and so $A_{1_0}=F=\langle \sigma_{\tau,a}\rangle$. Then $A_{1_0}$ has two orbits on $S_1$, that is, $U_1$ and $V_1$, and hence $\Gamma_{m,k,l}^{n}$ is not arc-transitive. To prove the half-arc-transitivity of $\Gamma_{m,k,l}^{n}$, we only need to show that $A$ is transitive on the vertex set and edge set of $\Gamma_{m,k,l}^{n}$. Note that $1_1\in U_1$ and $ (b^{m}a^{n})_1\in V_1$. To finish the proof, by Proposition~\ref{prop=Normalizer} it suffices to construct  a $\lambda\in \Aut(G)$ such that $\delta_{\lambda, b^{m}a^{n},1}\in I=\{\delta_{\lambda,x,y}\ |\ \lambda\in \Aut(G), S^\lambda=y^{-1}S^{-1}x\}$, that is, $S^{\lambda}=S^{-1}b^{m}a^{n}$, because $(1_0,1_1)^{\delta_{\lambda, b^{m}a^{n},1}}=((b^{m}a^{n})_1, 1_0)$.

Let $\mu=-(1+p^\gamma)^{m}-n(e-1)$ and $\nu=-(e-1)^{-1}\mu^{2}
-(e-1)^{-1}\mu$. Then $\mu+1+n(e-1)=0$ and hence
$\nu-\mu n
=-(e-1)^{-1}\mu^{2}-(e-1)^{-1}\mu-\mu n
=-(e-1)^{-1}\mu[\mu+1+n(e-1)]=0$ in  $\mz_{p^\gamma}$.
By Proposition~\ref{prop=element order}, $o(b^{m}a^{\nu-\mu n})=p^\b$. Denote by $m^{-1}$ the inverse of $m$ in $\mz_{p^\b}$. Then $(b^{m}a^{\nu-\mu n})^{m^{-1}}=ba^\ell$ for some $\ell$ in $\mz_{p^\a}$, and it is easy to check that $a^\mu$ and $(b^{m}a^{\nu-\mu n})^{m^{-1}}$ have the same relations as do $a$ and $b$. Define $\lambda$ as the automorphism of $G$ induced by $a\mapsto a^{\mu}$, $b\mapsto (b^{m}a^{\nu-\mu n})^{m^{-1}}$. Clearly, $(b^m)^\lambda=b^{m}a^{\nu-\mu n}$.

Note that $S=U\cup V$. Then $U^\lambda=\{a^{\eta\mu}\ |\ \eta\in T\}=\{a^\eta\ |\ \eta\in T\mu\}$ and $V^{-1}b^{m}a^{n}=\{(b^{m}a^{\eta})^{-1}b^{m}a^{n}\ |\ \eta\in T'\}=\{a^{-\eta+n}\ |\ \eta\in T'\}=\{a^\eta\ |\ \eta\in -T'+n\}$.
Recall that $T'=T[(1+p^\gamma)^{m}+n(e-1)]+n=-T\mu+n$. Then
$-T'+n=T\mu-n+n=T\mu$, and so $U^\lambda=V^{-1}b^{m}a^{n}$.
Similarly, $V^\lambda=\{(b^{m}a^{\eta})^\lambda\ |\ \eta\in T'\}
=\{b^{m}a^{\nu-\mu n}a^{\eta\mu}\ |\ \eta\in T'\}
=\{b^{m}a^{\eta}\ |\ \eta\in T^{'}\mu-\mu n+\nu\}$ and $U^{-1}b^{m}a^{n}=\{(a^\eta)^{-1}b^{m}a^{n}\ |\ \eta\in T\}=\{b^{m}a^{-\eta(1+p^\gamma)^{m}+n}\ |\ \eta\in T\}=\{b^{m}a^{\eta}\ |\ \eta\in -T(1+p^\gamma)^{m}+n\}$.
To prove $V^{\lambda}=U^{-1}b^{m}a^{n}$, we only need to show $T^{'}\mu-\mu n+\nu=-T(1+p^\gamma)^{m}+n$ in $\mz_{p^\a}$, which is equivalent to show that $T(1+p^\gamma)^{m}=T{\mu^2}-\nu+n$ because $T'=-T\mu+n$.

By Eq~(\ref{eq3}), $e^l(1+p^\gamma)^m=[(1+p^\gamma)^{m}-n(1-e)]^{2}=\mu^{2}$,
and by Lemma~\ref{lem=$e-1$ reversible and $Tx=T$}, $T=Te^l+(e-1)^{-1}(e^l-1)$. It follows
$T(1+p^\gamma)^m=Te^l(1+p^\gamma)^m+(e-1)^{-1}(e^l-1)(1+p^\gamma)^m=
T\mu^{2}+(e-1)^{-1}[\mu^{2}-(1+p^\gamma)^m]$. Note that $-\nu+n=(e-1)^{-1}\mu^{2}+(e-1)^{-1}\mu+n=(e-1)^{-1}[\mu^{2}+\mu+n(e-1)]
=(e-1)^{-1}[\mu^{2}-(1+p^\gamma)^m]$. Then $T(1+p^\gamma)^{m}=T{\mu^2}-\nu+n$, as required. Thus,  $V^{\lambda}=U^{-1}b^{m}a^{n}$ and hence $S^{\lambda}=S^{-1}b^{m}a^{n}$. This completes the proof. \hfill\qed

\medskip
\f {\bf Acknowledgements:} This work was supported by the National Natural Science Foundation of China (11571035) and by the 111 Project of China (B16002).

\end{document}